\documentclass[12pt]{article}
\usepackage{amsmath} 
\usepackage{graphicx}
\usepackage[all]{xy}
\usepackage{latexcad}
\input diagxy
\usepackage[usenames]{color}

\newtheorem{definition}{Definition}[section]

\newtheorem{prop}{Proposition}[section]
\newtheorem{theorem}{Theorem}[section]
\newtheorem{rem}{Remark}[section]

\newtheorem{thm}{Theorem}[section]
\newtheorem{cor}{Corollary}[section]

\newtheorem{dfn}{Definition}[section]

\font\bbfont=msbm10 at 12 pt

\def\bbR{\mbox{\bbfont R}}

\def\mathopdef#1{\expandafter\def\csname#1\endcsname{\mathop{\rm#1}\nolimits}}
\def\mathoplsdef#1{\expandafter\def\csname#1\endcsname{\mathop{\rm#1}}}
\def\mathbfdef#1{\expandafter\def\csname#1\endcsname{{\rm\bf#1}}}
\def\mathrmdef#1{\expandafter\def\csname#1\endcsname{{\rm#1}}}

\newcommand{\tn}{\otimes}

\newcommand{\ol}{\overline}               

\newcommand{\pf}{\vspace{\baselineskip}\noindent{\bf Proof. }}
\newcommand{\qed}{\hfill\rule{4pt}{8pt}\par\vspace{\baselineskip}}

\newcommand{\bC}{{\bf C}}

\newcommand{\eps}{{\epsilon}}

\newcommand{\tw}{{\tau}}

\newcommand{\set}{\mathbf{Set}}
\newcommand{\rel}{\mathbf{Rel}}

\newcommand{\gra}{\mathbf{Graph}}
\newcommand{\spa}{\mathbf{Span}}
\newcommand{\spg}{{\mathbf{Span}(\gra)}}
\newcommand{\cspg}{{\mathbf{Cospan}(\gra)}}
\newcommand{\csp}{{\mathbf {Cospan}}}
\newcommand{\gp}{\mathbf{Group}}
\newcommand{\tr}{\mathbf{TRel}}
\newcommand{\tsp}{\mathbf{TSpan}}
\newcommand{\tang}{\mathbf{Tangle}}

\newcommand{\tcd}{\mathbf{TCircD}}

\newcommand{\gpop}{{\mathbf {Group^{op}}}}
\title{Tangled circuits}
\author{R. Rosebrugh\\
Department of Mathematics and Computer Science\\
Mount Allison University\\
Sackville, N. B. E4L 1E6 Canada\\
{\tt rrosebrugh@mta.ca}\\
\and
N. Sabadini\\
Dipartimento di Informatica e Comunicazione\\
Universit\`a dell' Insubria\\
via Carloni, 78, Como, Italy\\
{\tt nicoletta.sabadini@uninsubria.it}
\and
R. F. C. Walters\\
Dipartimento di Informatica e Comunicazione\\
Universit\`a dell' Insubria\\
via Carloni, 78, Como, Italy\\
{\tt robert.walters@uninsubria.it}}

\begin{document}
\maketitle
\newpage

\section{Introduction}\label{sec-intro}

 The theme of the paper is the use of \emph{commutative Frobenius algebras in braided  
strict monoidal categories} in the study of varieties of circuits and communicating systems which occur in Computer Science, including circuits in which the wires are tangled. We indicate also some possible novel geometric interest in such algebras.

The contribution of the paper is the introduction and application of several new such categories, and appropriate functors between them. The authors and collaborators have previously studied similar systems using symmetric  monoidal categories (\cite{KSW97b,KSW00a, KSW00b, KSW02,RSW04,RSW05,RSW08,dFSWcsg}), with separable algebras instead of Frobenius algebras. These earlier works did not take into consideration any tangling of the wires. Further we will see in section \ref{sec-RLC} the importance of considering Frobenius algebras rather than the more special separable algebras even in the symmetric monoidal case (no tangling).

\subsection{Tangled circuit diagrams}\label{subsec-tangled-circuits}

We propose  a definition for a category of \emph{tangled circuit diagrams}, in which it is possible to distinguish, for example, the first and second of the following circuit diagrams, while the second and third are equal.

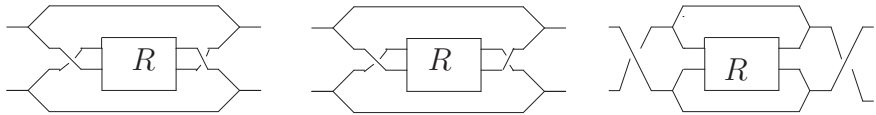
\begin{center}\label{fig-different-circuits}
{\tt\setlength{\unitlength}{2.0pt}
\begin{picture}(172,28)
\thinlines
\drawpath{4.13}{20.02}{8.13}{20.02}
\drawpath{8.13}{20.02}{12.13}{24.02}
\drawpath{12.13}{24.02}{44.13}{24.02}
\drawpath{44.13}{24.02}{48.13}{20.02}
\drawpath{48.13}{20.02}{52.13}{20.02}
\drawpath{8.13}{20.02}{12.13}{16.02}
\drawpath{12.13}{16.02}{14.13}{16.02}
\drawpath{14.13}{16.02}{18.13}{12.02}
\drawpath{18.13}{12.02}{22.13}{12.02}
\drawpath{22.13}{12.02}{22.13}{16.02}
\drawpath{22.03}{15.97}{22.03}{17.97}
\drawpath{36.13}{16.02}{36.13}{8.02}
\drawpath{22.13}{18.02}{36.13}{18.02}
\drawpath{22.13}{8.02}{22.13}{12.02}
\drawpath{36.13}{18.02}{36.13}{16.02}
\drawpath{22.13}{16.02}{18.13}{16.02}
\drawpath{18.13}{16.02}{18.13}{16.02}
\drawpath{18.29}{15.77}{17.06}{14.2}
\drawpath{14.07}{11.76}{15.84}{13.15}
\drawpath{14.13}{12.02}{12.13}{12.02}
\drawpath{12.13}{12.02}{8.13}{8.02}
\drawpath{8.13}{8.02}{4.13}{8.02}
\drawpath{4.13}{8.02}{6.13}{8.02}
\drawpath{6.13}{8.02}{8.13}{8.02}
\drawpath{8.13}{8.02}{12.13}{4.02}
\drawpath{12.13}{4.02}{44.13}{4.02}
\drawpath{44.13}{4.02}{48.13}{8.02}
\drawpath{48.13}{8.02}{52.13}{8.02}
\drawpath{52.13}{8.02}{48.13}{8.02}
\drawpath{48.13}{8.02}{44.13}{12.02}
\drawpath{48.13}{20.02}{44.13}{16.02}
\drawpath{44.13}{16.02}{42.13}{16.02}
\drawpath{44.13}{12.02}{42.13}{12.02}
\drawpath{40.13}{16.02}{42.13}{12.02}
\drawpath{40.13}{16.02}{36.13}{16.02}
\drawpath{36.13}{12.02}{40.13}{12.02}
\drawpath{22.13}{8.02}{36.13}{8.02}
\drawpath{40.79}{12.79}{39.91}{11.76}
\drawpath{42.38}{15.77}{41.86}{14.54}
\drawpath{61.79}{19.84}{65.79}{19.84}
\drawpath{65.79}{19.84}{69.79}{23.84}
\drawpath{69.79}{23.84}{101.79}{23.84}
\drawpath{101.79}{23.84}{105.79}{19.84}
\drawpath{105.79}{19.84}{109.79}{19.84}
\drawpath{65.79}{19.84}{69.79}{15.84}
\drawpath{69.79}{15.84}{71.79}{15.84}
\drawpath{71.79}{15.84}{75.79}{11.84}
\drawpath{75.79}{11.84}{79.79}{11.84}
\drawpath{79.79}{11.84}{79.79}{15.84}
\drawpath{79.79}{15.84}{79.79}{17.84}
\drawpath{93.79}{15.84}{93.79}{7.84}
\drawpath{79.79}{17.84}{93.79}{17.84}
\drawpath{79.79}{7.84}{79.79}{11.84}
\drawpath{93.79}{17.84}{93.79}{15.84}
\drawpath{79.79}{15.84}{75.79}{15.84}
\drawpath{75.79}{15.84}{75.79}{15.84}
\drawpath{75.94}{15.59}{74.72}{14.02}
\drawpath{71.73}{11.58}{73.49}{12.97}
\drawpath{71.79}{11.84}{69.79}{11.84}
\drawpath{69.79}{11.84}{65.79}{7.84}
\drawpath{65.79}{7.84}{61.79}{7.84}
\drawpath{61.79}{7.84}{63.79}{7.84}
\drawpath{63.79}{7.84}{65.79}{7.84}
\drawpath{65.79}{7.84}{69.79}{3.84}
\drawpath{69.79}{3.84}{101.79}{3.84}
\drawpath{101.79}{3.84}{105.79}{7.84}
\drawpath{105.79}{7.84}{109.79}{7.84}
\drawpath{109.79}{7.84}{105.79}{7.84}
\drawpath{105.79}{7.84}{101.79}{11.84}
\drawpath{105.79}{19.84}{101.79}{15.84}
\drawpath{101.79}{15.84}{99.79}{15.84}
\drawpath{101.79}{11.84}{99.79}{11.84}
\drawpath{100.0}{16.0}{98.0}{12.0}
\drawpath{97.79}{15.84}{93.79}{15.84}
\drawpath{93.79}{11.84}{97.79}{11.84}
\drawpath{79.79}{7.84}{93.79}{7.84}
\drawpath{97.83}{15.75}{98.54}{14.7}
\drawpath{99.25}{12.79}{99.94}{11.75}
\drawpath{132.13}{12.02}{136.13}{12.02}
\drawpath{136.13}{12.02}{136.13}{16.02}
\drawpath{136.03}{15.97}{136.03}{17.97}
\drawpath{150.13}{16.02}{150.13}{8.02}
\drawpath{136.13}{18.02}{150.13}{18.02}
\drawpath{136.13}{8.02}{136.13}{12.02}
\drawpath{150.13}{18.02}{150.13}{16.02}
\drawpath{136.13}{16.02}{132.13}{16.02}
\drawpath{132.13}{16.02}{132.13}{16.02}
\drawpath{154.13}{16.02}{150.13}{16.02}
\drawpath{150.13}{12.02}{154.13}{12.02}
\drawpath{136.13}{8.02}{150.13}{8.02}
\drawpath{132.0}{16.0}{130.0}{20.0}
\drawpath{130.0}{20.0}{132.0}{24.0}
\drawpath{132.0}{22.0}{132.0}{22.0}
\drawpath{132.0}{24.0}{154.0}{24.0}
\drawpath{154.0}{24.0}{156.0}{20.0}
\drawpath{156.0}{20.0}{154.0}{16.0}
\drawpath{132.0}{12.0}{130.0}{8.0}
\drawpath{130.0}{8.0}{132.0}{4.0}
\drawpath{132.0}{4.0}{154.0}{4.0}
\drawpath{154.0}{4.0}{156.0}{8.0}
\drawpath{156.0}{8.0}{154.0}{12.0}
\drawpath{130.0}{8.0}{126.0}{8.0}
\drawpath{126.0}{8.0}{120.0}{20.0}
\drawpath{120.0}{20.0}{118.0}{20.0}
\drawpath{118.0}{8.0}{120.0}{8.0}
\drawpath{120.0}{8.0}{122.0}{14.0}
\drawpath{124.0}{16.0}{126.0}{20.0}
\drawpath{126.0}{20.0}{130.0}{20.0}
\drawpath{156.0}{20.0}{160.0}{20.0}
\drawpath{156.0}{8.0}{160.0}{8.0}
\drawpath{160.0}{8.0}{166.0}{20.0}
\drawpath{166.0}{20.0}{168.0}{20.0}
\drawpath{160.0}{20.0}{162.0}{14.0}
\drawpath{164.0}{12.0}{166.0}{6.0}
\drawpath{166.0}{6.0}{168.0}{6.0}
\drawcenteredtext{30.0}{14.0}{$R$}
\drawcenteredtext{86.0}{14.0}{$R$}
\drawcenteredtext{142.0}{12.0}{$R$}
\end{picture}
}
\end{center}

The notion of tangled circuit diagram is parametrized by a multigraph (or tensor scheme) of components (such as the component $R$ in the example above). Given such a multigraph $M$, a tangled circuit diagram (or more briefly, a circuit diagram) is an arrow in the free braided strict monoidal category on $M$ in which objects of the multigraph $M$ are equipped with symmetric Frobenius algebra structures; we denote this category by $\tcd_M$. The objects of the multigraph $M$ may be thought of as \emph{types of wires}. Given any object $A$ of $M$ it is straightforward to see that there is an appropriate functor from Freyd, Yetter's category $\tang$ (\cite{FY89}) to $\tcd_M$ since a symmetric Frobenius structure on $A$ induces a tangle algebra structure on $A$. As a result any invariants of tangled circuit diagrams provide also invariants for tangles and knots. We conjecture that such functors $\tang \to \tcd_M$ are faithful.  We also conjecture that there is a topological description of $\tcd_M$ related to Freyd, Yetter's description of $\tang$ and to cobordisms. 

\subsection{Relations}
The category $\rel$ whose objects are sets, and whose arrows are relations is symmetric monoidal with the tensor of sets being the cartesian product, and each object has a symmetric Frobenius (even separable) algebra structure provided by the diagonal functions and their reverse relations.
In fact this was the motivating example for the introduction in  \cite{CW87} of the Frobenius equations  (equivalent axioms had been given earlier by Lawvere in \cite{L67}).
We describe here a modification of $\rel$ which we call  $\tr_G$, which depends on a group $G$, and which is braided rather than symmetric. We further describe a commutative Frobenius algebra in $\tr_G$ which hence yields a representation of $\tcd_M$, and this representation enables us, for example, to distinguish the two different circuits above. We discuss distinguishing closed circuits, a problem analogous to classifying knots, using $\tr_G$.

\subsection{Spans and cospans}
The principal category we have using in the earlier work on circuits and communicating-parallel algebras of processes has been the category \newline $\spg$ of spans of graphs (and for sequential systems \newline $\cspg$).  Already in the original paper \cite{KSW97b} the separable algebra structure on each object played a crucial role. The relation between another model of circuits, namely Mealy automata and $\spg$ was discussed in \cite{KSW00a}. One of the motivations of the present work is to produce an semantic algebra in which the twisting of wires is also (at least partially) expressible. To this end we introduce first a simple braided modification $\tsp_G$ of $\spa (\set)$, depending on a group $G$, with a commutative Frobenius algebra. It is clear that a similar construction $\tsp_G(\bC)$ could be made for a group object $G$ in a category $\bC$ with limits in the place of $\set$.

Again, there is a representation of $\tang$ (via a representation of $\tcd$) which takes a tangle to the span of colourings of the tangle (introduced by John Armstrong in \cite{accn}). Applied to knots the set of colourings is one of the simplest invariants for distinguishing knots (as a first example it allows one to show that a trefoil is not an unknot). The extended notion of colourings of tangled circuit diagrams gives further aid in distinguishing circuit diagrams.

The category $\gpop$, the dual of the category of groups has finite limits. Further $F$, the free group on one generator is a group object in $\gpop$. The category $\tsp_F (\gpop)$ is braided monoidal with $F$ equipped with a commutative Frobenius structure. The induced representation $$\tang \to \tsp_F (\gpop)$$ associates the cospan of groups introduced by John Armstrong in \cite{aekt}  to a tangle, and the knot group to a knot.

\subsection{Linear analogue circuits}
This example comes from the paper \cite{KSW00a} where it is discussed in detail. However 
the Frobenius algebra structure was not noticed in that paper. The category is analogous to $\tsp_G (\gra)$ where the group $G$ is the real numbers under addition. The Frobenius algebra structure arises from the Kirchhoff law for currents. Since the group is abelian there is no information about the tangling of wires. We describe, as an example, circuits composed of resistors, capacitors and inductors.  

\subsection{Remaining questions}
Proving that two expressions in $\tcd$ yield different circuit seems to be a difficult question even in apparently simple cases some of which we note below. If as we suspect knots are faithfully represented in $\tcd$ this is not surprising, though for knots there are known though non-trivial algorithms.

\subsection{References}
There is a huge literature now relating monoidal categories and geometry beginning with \cite{p71,KL80,jsbmc}. We mention just two further
items of an expository nature useful to reading this paper (apart from our own work mentioned above): the first \cite{S11} is a survey for computer scientists and others which discusses many additional structures but strangely not Frobenius algebras, and ignores our work on separable algebras; the second \cite{K04} is an introductory book on the relation between Frobenius algebras and 2-dimensional cobordism. 

\subsection{Acknowledgement} 
We would like to thank Aurelio Carboni for helpful suggestions.
\section{Braided monoidal categories and Frobenius algebras}\label{sec-braided-monoidal}
We review immediately the notions fundamental for the paper.
\subsection{Braided monoidal categories}
\begin{definition}
 A \emph{braided strict monoidal category} (\cite{jsbmc}) is a category $\bC$ with a functor, called tensor,
$\tn : \bC \times \bC \to \bC$ and a ``unit'' object $I$
together with a natural family of isomorphisms $\tw_{A,B} : A \tn B \to B \tn A$ called twist
satisfying
\begin{itemize}
\item[1)] $\tn$ is associative and unitary on objects and arrows,
\item[2)] the following diagrams commute for objects $A, B, C$: 
$$
\bfig
\place(-400,200)[B1:]
\Vtriangle/>`>`<-/<600,400>[A\tn B\tn C`B\tn C\tn A`B\tn A\tn C;\tw`\tw\tn 1`1\tn\tw]
\efig
$$
and
$$
\bfig
\place(-400,200)[B2:]
\Vtriangle/>`>`<-/<600,400>[A\tn B\tn C`C\tn A\tn B`A\tn C\tn B;\tw`1\tn\tw`\tw\tn 1]
\efig
$$
\end{itemize}
\end{definition}

Among the consequences of the definition is the Yang-Baxter equation which  reads:
$$(1\tn\tw)(\tw\tn 1)(1\tn\tw) = (\tw\tn 1)(1\tn\tw)(\tw\tn 1) : A\tn B\tn C \to C\tn B\tn A$$
A compact and comprehensible formulation of such properties
is provided by circuit or``wire'' diagrams like the following.
Composition is read from left to right and $\tn$ is vertical juxtaposition. The twist
is expressed by the ``positive crossing'' (top wire over bottom) and its inverse by the negative crossing.
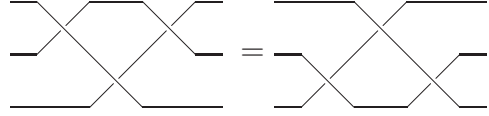
\begin{center}
\begin{picture}(200,70)
\put(0,0){\line(1,0){30}}
\put(30,0){\line(1,1){9}}
\put(41,11){\line(1,1){18}}
\put(61,31){\line(1,1){9}}
\put(70,40){\line(1,0){10}}
 
\put(0,20){\line(1,0){10}}
\put(10,20){\line(1,1){9}}
\put(21,31){\line(1,1){9}}
\put(30,40){\line(1,0){20}}
\put(50,40){\line(1,-1){20}}
\put(70,20){\line(1,0){10}}
 
\put(0,40){\line(1,0){10}}
\put(10,40){\line(1,-1){40}}
\put(50,0){\line(1,0){30}}
 
\put(87,18){=}
 
\put(100,0){\line(1,0){10}}
\put(110,0){\line(1,1){9}}
\put(121,11){\line(1,1){18}}
\put(141,31){\line(1,1){9}}
\put(150,40){\line(1,0){30}}
 
\put(100,20){\line(1,0){10}}
\put(110,20){\line(1,-1){20}}
\put(130,0){\line(1,0){20}}
\put(150,0){\line(1,1){9}}
\put(161,11){\line(1,1){9}}
\put(170,20){\line(1,0){10}}
 
\put(100,40){\line(1,0){30}}
\put(130,40){\line(1,-1){40}}
\put(170,0){\line(1,0){10}}
 
\end{picture}
\end{center}


Another consequence of the axioms above is that $\tau_{A,I}=\tau_{I,A}=1_A:A\to A$.

The naturality of the twist $\tau$ leads to the following kind of equality of diagrams:
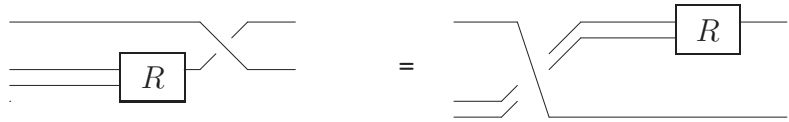
\begin{center}
{\tt\setlength{\unitlength}{3.0pt}
\begin{picture}(100,25)      
\thinlines
\drawframebox{22.0}{9.0}{8.0}{6.0}{$R$}
\drawpath{18.0}{10.0}{10.0}{10.0}
\drawpath{18.0}{8.0}{10.0}{8.0}
\drawpath{10.0}{10.0}{4.0}{10.0}
\drawpath{4.0}{6.0}{4.0}{6.0}
\drawpath{10.0}{8.0}{4.0}{8.0}
\drawpath{4.0}{16.0}{28.0}{16.0}
\drawpath{26.0}{10.0}{28.0}{10.0}
\drawpath{96.0}{16.0}{102.0}{16.0}
\drawpath{28.0}{10.0}{30.0}{12.0}
\drawpath{28.0}{16.0}{34.0}{10.0}
\drawpath{24.0}{72.0}{24.0}{72.0}
\drawcenteredtext{54.0}{10.0}{=}
\drawpath{60.0}{16.0}{68.0}{16.0}
\drawpath{68.0}{16.0}{72.0}{4.0}
\drawpath{72.0}{4.0}{102.0}{4.0}
\drawpath{60.0}{6.0}{66.0}{6.0}
\drawpath{66.0}{6.0}{68.0}{8.0}
\drawpath{60.0}{4.0}{66.0}{4.0}
\drawpath{66.0}{4.0}{68.0}{6.0}
\drawpath{72.0}{10.0}{76.0}{14.0}
\drawpath{72.0}{12.0}{76.0}{16.0}
\drawpath{76.0}{16.0}{88.0}{16.0}
\drawpath{76.0}{14.0}{88.0}{14.0}
\drawframebox{92.0}{15.0}{8.0}{6.0}{$R$}
\drawpath{24.0}{72.0}{24.0}{72.0}
\drawpath{32.0}{14.0}{34.0}{16.0}
\drawpath{34.0}{16.0}{40.0}{16.0}
\drawpath{34.0}{10.0}{40.0}{10.0}
\end{picture}}
\end{center}   
In the case when the codomain of the component is $I$ naturality is drawn, for example, as:

\begin{center}
{\tt\setlength{\unitlength}{3.0pt} 
\begin{picture}(106,22)
\thinlines
\drawframebox{22.0}{9.0}{8.0}{6.0}{$R$}
\drawpath{18.0}{10.0}{10.0}{10.0}
\drawpath{18.0}{8.0}{10.0}{8.0}
\drawpath{10.0}{10.0}{4.0}{10.0}
\drawpath{4.0}{6.0}{4.0}{6.0}
\drawpath{10.0}{8.0}{4.0}{8.0}
\drawpath{4.0}{16.0}{28.0}{16.0}
\drawpath{28.0}{16.0}{38.0}{16.0}
\drawcenteredtext{54.0}{10.0}{=}
\drawpath{60.0}{16.0}{68.0}{16.0}
\drawpath{68.0}{16.0}{72.0}{4.0}
\drawpath{72.0}{4.0}{102.0}{4.0}
\drawpath{60.0}{6.0}{66.0}{6.0}
\drawpath{66.0}{6.0}{68.0}{8.0}
\drawpath{60.0}{4.0}{66.0}{4.0}
\drawpath{66.0}{4.0}{68.0}{6.0}
\drawpath{72.0}{10.0}{76.0}{14.0}
\drawpath{72.0}{12.0}{76.0}{16.0}
\drawpath{76.0}{16.0}{88.0}{16.0}
\drawpath{76.0}{14.0}{88.0}{14.0}
\drawframebox{92.0}{15.0}{8.0}{6.0}{$R$}
\end{picture}}
\end{center}

\subsubsection{Frobenius algebras}
\begin{dfn}\label{def-Frobenius-algebra}
A \emph{commutative Frobenius algebra} in a braided monoidal category consists of an object $G$
and four arrows $\nabla :G\otimes G\to G$, $\Delta :G\to G\otimes G$, $n:I\to G$ and $e:G\to I$
making $(G,\nabla,e)$ a monoid, $(G,\Delta,n)$ a comonoid and satisfying the equations
\begin{align}
(1_{G}\tn\nabla)(\Delta\tn 1_{G}) &  =\Delta\nabla = (\nabla\tn 1_{G})(1_{G}\tn\Delta) : G\tn G\to G\tn G \tag{D}\\
 \nabla\tw & = \nabla : G\tn G\to G \notag\\
 \tw\Delta & =\Delta : G\to G\tn G \notag
\end{align}
\end{dfn}

\begin{dfn}
A multigraph $M$ consists of two sets $M_0$ (vertices or wires) and $M_1$ (edges or components) and two functions $dom:M_1\to M_0^*$ and $cod:M_1\to M_0^*$ where $M_0^*$ is the free monoid on $M_0$.
\end{dfn}

\begin{dfn}\label{def-TangCircDiag}
Given a multigraph $M$ the free braided strict monoidal category in which the objects of $M$ are equipped with commutative Frobenius algebra structures is called $\tcd_M$. Its arrows are called tangled circuit diagrams, or more briefly circuit diagrams. In the case that $M$ has one vertex and no arrows we will denote $\tcd_M$ simply as $\tcd$.
\end{dfn}

\subsubsection{Tangle algebras}
\begin{definition}\label{def-tangle-algebra}
An object $X$ in a braided strict monoidal category (with twist $\tau$)is called a \emph{tangle algebra} when it is equipped with
arrows $\eta : I \to X \tn X$ and $\eps : X \tn X \to I$ that satisfy the equations (where we write 
$1$ for all identities):
\begin{itemize}
\item[(i)] $(\eps\tn 1) (1\tn \eta) = 1 = (1\tn \eps)(\eta\tn 1) $
\item[(ii)]  $\eps \tw = \eps$ and $\tw\eta  = \eta$.
\end{itemize}

Axiom (i) says that $X$ is a \emph{self-dual object}.
The reader can translate these into wire diagrams.
An example is the wire diagram for (i):

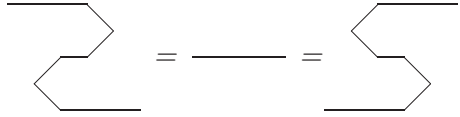
\begin{center}
\begin{picture}(200,70)(-40,0)
\put(-10,60){\line(1,0){30}}
\put(0,30){\line(1,-1){10}}
\put(0,30){\line(1,1){10}}
\put(20,40){\line(1,1){10}}
\put(10,40){\line(1,0){10}}
\put(20,60){\line(1,-1){10}}
\put(10,20){\line(1,0){30}}
\put(45,37){=}
\put(60,40){\line(1,0){35}}
\put(100,37){=}
\put(130,60){\line(1,0){30}}
\put(120,50){\line(1,-1){10}}
\put(120,50){\line(1,1){10}}
\put(130,40){\line(1,0){10}}
\put(140,40){\line(1,-1){10}}
\put(140,20){\line(1,1){10}}
\put(110,20){\line(1,0){30}}
\end{picture}
\end{center}
\end{definition}

\begin{thm} If $G$ is a commutative Frobenius algebra in a braided monoidal category, then the arrows
$\eps = e\nabla$, $\eta = \Delta n$, $\tw$
satisfy the axioms of the generating object of the category of tangles, and so $G$ is a tangle algebra.
\end{thm}

\pf
Let $G$ be a commutative Frobenius algebra in a braided monoidal category. It is straightforward to give
algebraic proofs for the tangle algebra axioms, but we remind the reader that these
can be more easily found using wire diagrams. 

To see that $(\eps\tn 1)(1\tn \eta)=1 $ notice that 
$$(e\tn 1)(\nabla\tn 1)(1\tn\Delta)(1\tn n) = (e\tn 1)\Delta\nabla(1\tn n)=1\cdot 1=1.$$
\qed

\begin{definition}\label{def-tang} (Freyd-Yetter)
The category $\tang$ is the free strict monoidal category generated by one object $X$,
equipped with a tangle algebra structure.
\end{definition}

The category $\tang$ has a geometric description \cite{Y2001} consonant with its name. In that description the arrows from $I$ to $I$ are knots and links.

\begin{cor}\label{cor-tang-to-trel}
Given an object $A$ of multigraph $M$ there is a unique braided strict monoidal functor $\tang \to \tcd_M$ taking the generating object to
$A$ and the structure maps of $\tang$ to the corresponding structure maps of $A$ in $\tcd_M$.
\end{cor}

\subsubsection{Example equations}\label{subsubsec-single-twist}

We now give some examples of equations between circuit diagrams.

\begin{prop}
If $R:X\times Y\to I$ is an arrow in the multigraph $M$ then 
$$R\tau_{Y,X}^{-1}=(\eps_X)(1_X\tn R\tn 1_X)(\eta\tn 1_Y\tn 1_X)=R\tau_{X,Y}$$.
\end{prop}

\pf

First a picture of the equations:

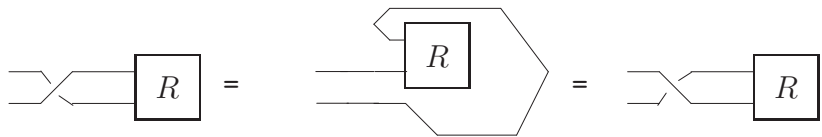
\begin{center}
{\tt\setlength{\unitlength}{3.0pt} 
\begin{picture}(110,24)
\thinlines
\drawpath{12.0}{12.0}{20.0}{12.0}
\drawpath{12.0}{8.0}{20.0}{8.0}
\drawframebox{24.0}{10.0}{8.0}{8.0}{$R$}
\drawcenteredtext{32.0}{10.0}{=}
\drawpath{42.69}{12.02}{46.69}{12.02}
\drawpath{45.79}{12.02}{50.02}{12.02}
\drawpath{50.0}{8.0}{54.0}{8.0}
\drawpath{42.81}{8.01}{46.81}{8.01}
\drawpath{46.6}{8.01}{50.6}{8.01}
\drawpath{8.0}{8.0}{12.0}{12.0}
\drawpath{50.16}{12.02}{54.16}{12.02}
\drawpath{54.0}{8.0}{58.0}{4.0}
\drawpath{58.0}{4.0}{68.0}{4.0}
\drawframebox{58.0}{14.0}{8.0}{8.0}{$R$}
\drawpath{54.0}{16.0}{52.0}{16.0}
\drawpath{52.0}{16.0}{50.0}{18.0}
\drawpath{50.0}{18.0}{52.0}{20.0}
\drawpath{52.0}{20.0}{66.0}{20.0}
\drawpath{66.0}{20.0}{72.0}{12.0}
\drawpath{72.0}{12.0}{68.0}{4.0}
\drawpath{8.0}{8.0}{4.0}{8.0}
\drawpath{8.0}{12.0}{4.0}{12.0}
\drawpath{8.01}{12.1}{9.25}{10.53}
\drawpath{10.3}{8.97}{12.06}{7.75}
\drawpath{90.0}{12.0}{98.0}{12.0}
\drawpath{90.0}{8.0}{98.0}{8.0}
\drawframebox{102.0}{10.0}{8.0}{8.0}{$R$}
\drawpath{86.0}{12.0}{90.0}{8.0}
\drawpath{86.0}{8.0}{82.0}{8.0}
\drawpath{86.0}{12.0}{82.0}{12.0}
\drawpath{85.8}{7.75}{86.87}{9.31}
\drawpath{90.04}{11.93}{88.44}{10.88}
\drawcenteredtext{76.0}{10.0}{=}
\end{picture}
}

\end{center}

It is clearly sufficient to prove the first equation.

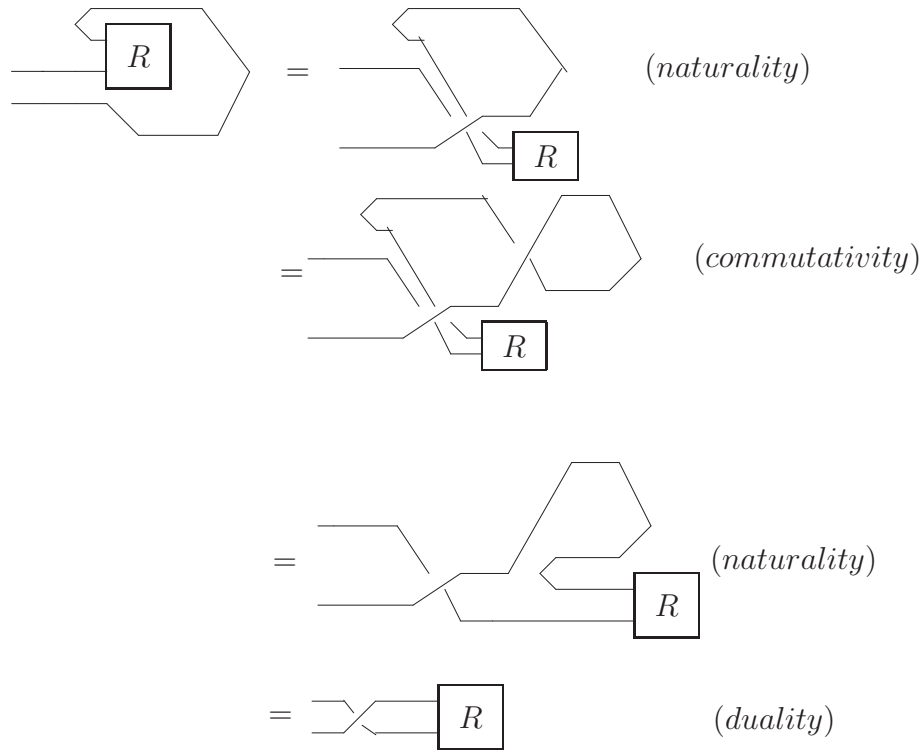
\begin{center}
{\tt\setlength{\unitlength}{3.0pt}
\begin{picture}(117,53)
\thinlines
\drawpath{54.0}{8.0}{42.0}{8.0}
\drawcenteredtext{105.0}{18.0}{$(commutativity)$}
\drawpath{8.5}{41.65}{12.65}{41.65}
\drawpath{4.63}{41.65}{8.65}{41.65}
\drawpath{12.65}{37.58}{16.65}{37.58}
\drawpath{4.65}{37.58}{8.65}{37.58}
\drawpath{8.65}{37.58}{12.65}{37.58}
\drawpath{12.65}{41.65}{16.65}{41.65}
\drawpath{16.65}{37.58}{20.65}{33.58}
\drawpath{20.65}{33.58}{30.65}{33.58}
\drawframebox{20.65}{43.58}{7.98}{8.0}{$R$}
\drawpath{16.65}{45.58}{14.65}{45.58}
\drawpath{14.65}{45.58}{12.65}{47.58}
\drawpath{12.65}{47.58}{14.65}{49.58}
\drawpath{14.65}{49.58}{28.65}{49.58}
\drawpath{28.65}{49.58}{34.65}{41.58}
\drawpath{34.65}{41.58}{30.65}{33.58}
\drawcenteredtext{40.0}{16.0}{$=$}
\drawpath{74.0}{26.0}{66.0}{12.0}
\drawpath{64.0}{26.0}{68.0}{20.0}
\drawpath{70.0}{18.0}{72.0}{14.0}
\drawpath{74.0}{26.0}{80.0}{26.0}
\drawpath{80.0}{26.0}{84.0}{18.0}
\drawpath{84.0}{18.0}{80.0}{14.0}
\drawpath{80.0}{14.0}{72.0}{14.0}
\drawcenteredtext{40.75}{41.24}{$=$}
\drawpath{48.0}{42.0}{46.0}{42.0}
\drawpath{56.0}{46.0}{62.0}{36.0}
\drawpath{74.0}{42.0}{70.0}{36.0}
\drawpath{70.0}{36.0}{64.0}{36.0}
\drawpath{64.0}{36.0}{58.0}{32.0}
\drawpath{56.0}{42.0}{60.0}{36.0}
\drawpath{58.0}{32.0}{46.0}{32.0}
\drawpath{48.0}{42.0}{56.0}{42.0}
\drawpath{64.0}{34.0}{66.0}{32.0}
\drawpath{56.65}{45.58}{54.65}{45.58}
\drawpath{54.65}{45.58}{52.65}{47.58}
\drawpath{52.65}{47.58}{54.65}{49.58}
\drawpath{54.65}{49.58}{68.66}{49.58}
\drawpath{68.66}{49.58}{74.66}{41.58}
\drawpath{62.0}{34.0}{64.0}{30.0}
\drawpath{66.0}{32.0}{68.0}{32.0}
\drawpath{64.0}{30.0}{68.0}{30.0}
\drawframebox{72.0}{31.0}{8.0}{6.0}{$R$}
\drawcenteredtext{95.0}{42.0}{$(naturality)$}
\drawpath{44.0}{18.0}{42.0}{18.0}
\drawpath{52.0}{22.0}{58.0}{12.0}
\drawpath{66.0}{12.0}{60.0}{12.0}
\drawpath{60.0}{12.0}{54.0}{8.0}
\drawpath{52.0}{18.0}{56.0}{12.0}
\drawpath{44.0}{18.0}{52.0}{18.0}
\drawpath{60.0}{10.0}{62.0}{8.0}
\drawpath{52.65}{21.58}{50.65}{21.58}
\drawpath{50.65}{21.58}{48.65}{23.58}
\drawpath{48.65}{23.58}{50.65}{25.58}
\drawpath{50.65}{25.58}{64.66}{25.58}
\drawpath{58.0}{10.0}{60.0}{6.0}
\drawpath{62.0}{8.0}{64.0}{8.0}
\drawpath{60.0}{6.0}{64.0}{6.0}
\drawframebox{68.0}{7.0}{8.0}{6.0}{$R$}
\end{picture}
}
\end{center}

\begin{center}
{\tt\setlength{\unitlength}{3.0pt}
\begin{picture}(50,44)
\thinlines
\drawpath{17.04}{10.83}{25.04}{10.83}
\drawpath{17.04}{6.84}{25.04}{6.84}
\drawframebox{29.04}{8.84}{8.0}{8.0}{$R$}
\drawcenteredtext{5.6}{28.18}{$=$}
\drawpath{41.77}{40.95}{33.77}{26.95}
\drawpath{39.77}{28.95}{37.77}{26.95}
\drawpath{37.77}{26.95}{39.77}{24.95}
\drawpath{41.77}{40.95}{47.77}{40.95}
\drawpath{13.04}{6.84}{17.04}{10.84}
\drawpath{47.77}{40.95}{51.77}{32.95}
\drawpath{51.77}{32.95}{47.77}{28.95}
\drawpath{47.77}{28.95}{39.77}{28.95}
\drawpath{11.77}{32.95}{9.77}{32.95}
\drawpath{31.77}{20.95}{49.77}{20.95}
\drawpath{33.77}{26.95}{27.77}{26.95}
\drawpath{27.77}{26.95}{21.77}{22.95}
\drawpath{19.77}{32.95}{23.77}{26.95}
\drawpath{21.77}{22.95}{9.77}{22.95}
\drawpath{11.77}{32.95}{19.77}{32.95}
\drawpath{13.04}{6.84}{9.04}{6.84}
\drawpath{13.04}{10.84}{9.04}{10.84}
\drawpath{13.04}{10.95}{14.29}{9.38}
\drawpath{15.34}{7.81}{17.09}{6.59}
\drawpath{25.77}{24.95}{27.77}{20.95}
\drawpath{27.77}{20.95}{31.77}{20.95}
\drawpath{39.77}{24.95}{49.77}{24.95}
\drawframebox{53.77}{22.95}{8.0}{8.0}{$R$}
\drawcenteredtext{5.04}{8.84}{$=$}
\drawcenteredtext{69.76}{28.54}{$(naturality)$}
\drawcenteredtext{67.3}{8.33}{$(duality)$}
\end{picture}
}
\end{center}
\qed

\begin{prop}
If $R:I\to X\otimes X$ and  $S:X\otimes X \to I$ then $S\tau^{2n}R=SR$; that is $R$ and $S$ joined
by an even number of twists is equal to $R$ and $S$ joined directly.
\end{prop}

A couple more equations provable in $\tcd$;
 
\subsubsection{Example}
\begin{center}
{\tt\setlength{\unitlength}{2.0pt}
\begin{picture}(130,26)
\thinlines
\drawpath{4.0}{16.0}{14.0}{16.0}
\drawpath{14.0}{16.0}{22.0}{22.0}
\drawpath{22.0}{22.0}{40.0}{22.0}
\drawpath{40.0}{22.0}{50.0}{16.0}
\drawpath{50.0}{16.0}{62.0}{16.0}
\drawpath{14.0}{16.0}{22.0}{10.0}
\drawpath{22.0}{10.0}{28.0}{10.0}
\drawpath{28.0}{10.0}{36.0}{4.0}
\drawpath{36.0}{4.0}{62.0}{4.0}
\drawpath{4.0}{4.0}{22.0}{4.0}
\drawpath{22.0}{4.0}{28.0}{4.0}
\drawpath{28.0}{4.0}{30.0}{6.0}
\drawpath{32.0}{8.0}{34.0}{10.0}
\drawpath{34.0}{10.0}{40.0}{10.0}
\drawpath{40.0}{10.0}{50.0}{16.0}
\drawpath{88.0}{22.0}{96.0}{22.0}
\drawpath{96.0}{22.0}{104.0}{16.0}
\drawpath{104.0}{16.0}{96.0}{10.0}
\drawpath{96.0}{10.0}{88.0}{10.0}
\drawpath{104.0}{16.0}{112.0}{16.0}
\drawpath{112.0}{16.0}{118.0}{22.0}
\drawpath{118.0}{22.0}{126.0}{22.0}
\drawpath{112.0}{16.0}{118.0}{10.0}
\drawpath{118.0}{10.0}{126.0}{10.0}
\drawcenteredtext{76.0}{16.0}{=}
\end{picture}
}
\end{center}

\begin{rem}
The geometric intuition is that the wires are \emph{thick} and so can be deformed contracting segments. Notice however that it is not true in general in $\tcd$ that the separable axiom  $\nabla \Delta=1$ holds. That is, cycles cannot be contracted to a point.
\end{rem}

\subsubsection{Example}\label{subsubsec-two-utilities}
\begin{center}
{\tt\setlength{\unitlength}{3.0pt}
\begin{picture}(98,32)
\thinlines
\drawframebox{7.21}{22.87}{5.97}{5.92}{}
\drawframebox{7.0}{7.0}{6.0}{6.0}{}
\drawframebox{29.0}{23.0}{6.0}{6.0}{}
\drawframebox{29.0}{7.0}{6.0}{6.0}{}
\drawpath{10.0}{24.0}{26.0}{24.0}
\drawpath{10.0}{22.0}{14.0}{22.0}
\drawpath{26.0}{6.0}{22.0}{6.0}
\drawpath{26.0}{8.0}{22.0}{8.0}
\drawpath{26.0}{22.0}{22.0}{22.0}
\drawpath{10.0}{8.0}{14.0}{8.0}
\drawpath{10.0}{6.0}{26.0}{6.0}
\drawpath{14.0}{22.0}{22.0}{8.0}
\drawpath{14.25}{8.06}{17.24}{13.65}
\drawpath{22.33}{21.84}{19.0}{16.25}
\drawcenteredtext{42.0}{14.0}{=}
\drawcenteredtext{8.0}{22.0}{$H_1$}
\drawcenteredtext{6.8}{6.2}{$H_2$}
\drawcenteredtext{30.0}{22.0}{$U_1$}
\drawcenteredtext{30.0}{7.4}{$U_2$}
\drawframebox{61.21}{22.87}{5.98}{5.92}{}
\drawframebox{61.0}{7.0}{6.0}{6.0}{}
\drawframebox{83.0}{23.0}{6.0}{6.0}{}
\drawframebox{83.0}{7.0}{6.0}{6.0}{}
\drawpath{64.0}{24.0}{80.0}{24.0}
\drawpath{64.0}{22.0}{68.0}{22.0}
\drawpath{80.0}{6.0}{76.0}{6.0}
\drawpath{80.0}{8.0}{76.0}{8.0}
\drawpath{80.0}{22.0}{76.0}{22.0}
\drawpath{64.0}{8.0}{68.0}{8.0}
\drawpath{64.0}{6.0}{80.0}{6.0}
\drawpath{68.0}{22.0}{76.0}{8.0}
\drawcenteredtext{62.0}{22.0}{$H_1$}
\drawcenteredtext{60.8}{6.2}{$H_2$}
\drawcenteredtext{84.0}{22.0}{$U_1$}
\drawcenteredtext{84.0}{7.4}{$U_2$}
\drawpath{68.0}{8.0}{70.0}{12.0}
\drawpath{70.0}{12.0}{68.0}{16.0}
\drawpath{68.0}{16.0}{54.0}{16.0}
\drawpath{54.0}{16.0}{50.0}{22.0}
\drawpath{50.0}{22.0}{54.0}{28.0}
\drawpath{54.0}{28.0}{88.0}{28.0}
\drawpath{88.0}{28.0}{94.0}{22.0}
\drawpath{94.0}{22.0}{90.0}{16.0}
\drawpath{90.0}{16.0}{76.0}{16.0}
\drawpath{76.0}{16.0}{74.0}{20.0}
\drawpath{74.0}{20.0}{76.0}{22.0}
\end{picture}
}
\end{center}

\subsubsection{Example}\label{subsubsec-flash}

If $R:I\to X\tn X$ and $S:X\tn X\to I$ then $$S(\eps\tn\eps\tn 1\tn 1)(1\tn \tau^{-1}\tn\tau^{-1}\tn 1)(1\tn1\tn \eta\tn \eta)R=S\tau\tau R.$$ 

Diagrammatically:

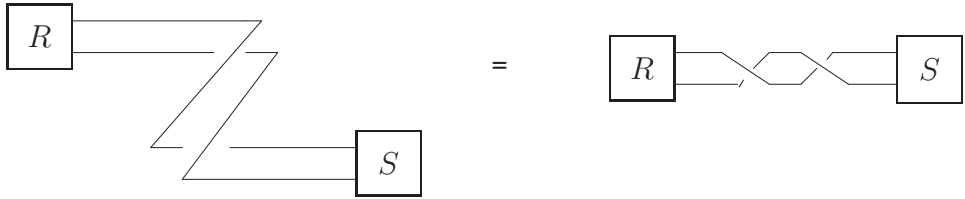
\begin{center}
{\tt\setlength{\unitlength}{3.0pt}
\begin{picture}(128,32)
\thinlines
\drawframebox{8.0}{24.0}{8.0}{8.0}{$R$}
\drawpath{12.0}{26.0}{34.0}{26.0}
\drawpath{34.0}{26.0}{36.0}{26.0}
\drawpath{36.0}{26.0}{22.0}{10.0}
\drawpath{22.0}{10.0}{26.0}{10.0}
\drawpath{12.0}{22.0}{30.0}{22.0}
\drawpath{34.0}{22.0}{36.0}{22.0}
\drawpath{36.0}{22.0}{38.0}{22.0}
\drawpath{38.0}{22.0}{26.0}{6.0}
\drawpath{26.0}{6.0}{48.0}{6.0}
\drawpath{32.0}{10.0}{48.0}{10.0}
\drawframebox{52.0}{8.0}{8.0}{8.0}{$S$}
\drawcenteredtext{66.0}{20.0}{=}
\drawframebox{84.0}{20.0}{8.0}{8.0}{$R$}
\drawpath{88.0}{22.0}{94.0}{22.0}
\drawpath{94.0}{22.0}{100.0}{18.0}
\drawpath{100.0}{18.0}{104.0}{18.0}
\drawpath{88.0}{18.0}{96.0}{18.0}
\drawpath{98.0}{20.0}{100.0}{22.0}
\drawpath{100.0}{22.0}{104.0}{22.0}
\drawpath{104.0}{22.0}{110.0}{18.0}
\drawpath{110.0}{18.0}{116.0}{18.0}
\drawpath{104.0}{18.0}{106.0}{20.0}
\drawpath{108.0}{22.0}{116.0}{22.0}
\drawframebox{120.22}{19.86}{8.07}{8.36}{$S$}
\drawpath{96.73}{18.52}{96.2}{17.65}
\drawpath{107.98}{22.0}{107.27}{21.13}
\end{picture}
}
\end{center}

\pf
We will give a diagrammatic proof. A more explicit picture of the left hand expression is

\begin{center}
\centerline{\tt\setlength{\unitlength}{2.5pt}
\begin{picture}(66,48)
\thinlines
\drawpath{42.0}{42.0}{14.0}{42.0}
\drawpath{28.0}{34.0}{14.0}{34.0}
\drawpath{34.0}{28.0}{42.0}{28.0}
\drawpath{42.0}{28.0}{44.0}{24.0}
\drawpath{44.0}{24.0}{42.0}{20.0}
\drawpath{42.0}{20.0}{34.0}{20.0}
\drawpath{34.0}{20.0}{28.0}{14.0}
\drawpath{28.0}{14.0}{24.0}{14.0}
\drawpath{24.0}{14.0}{22.0}{10.0}
\drawpath{22.0}{10.0}{24.0}{6.0}
\drawpath{24.0}{6.0}{52.0}{6.0}
\drawpath{42.0}{42.0}{44.0}{38.0}
\drawpath{44.0}{38.0}{42.0}{34.0}
\drawpath{42.0}{34.0}{34.0}{34.0}
\drawpath{34.0}{34.0}{28.0}{28.0}
\drawpath{28.0}{28.0}{24.0}{28.0}
\drawpath{24.0}{28.0}{22.0}{24.0}
\drawpath{22.0}{24.0}{24.0}{20.0}
\drawpath{24.0}{20.0}{28.0}{20.0}
\drawpath{34.0}{14.0}{52.0}{14.0}
\drawpath{28.0}{34.0}{30.0}{32.0}
\drawpath{32.0}{30.0}{34.0}{28.0}
\drawpath{28.0}{20.0}{30.0}{18.0}
\drawpath{32.0}{16.0}{34.0}{14.0}
\drawframebox{9.0}{38.0}{10.0}{12.0}{$R$}
\drawframebox{57.0}{10.0}{10.0}{12.0}{$S$}
\end{picture}
}
\end{center}

By naturality this is equal to

\begin{center}
\centerline{\tt\setlength{\unitlength}{2.0pt}
\begin{picture}(98,48)
\thinlines
\drawpath{74.0}{42.0}{46.0}{42.0}
\drawpath{60.0}{34.0}{46.0}{34.0}
\drawpath{66.0}{28.0}{74.0}{28.0}
\drawpath{74.0}{28.0}{76.0}{24.0}
\drawpath{76.0}{24.0}{74.0}{20.0}
\drawpath{74.0}{20.0}{66.0}{20.0}
\drawpath{66.0}{20.0}{60.0}{14.0}
\drawpath{60.0}{14.0}{56.0}{14.0}
\drawpath{56.0}{14.0}{54.0}{10.0}
\drawpath{54.0}{10.0}{56.0}{6.0}
\drawpath{56.0}{6.0}{84.0}{6.0}
\drawpath{74.0}{42.0}{76.0}{38.0}
\drawpath{76.0}{38.0}{74.0}{34.0}
\drawpath{74.0}{34.0}{66.0}{34.0}
\drawpath{66.0}{34.0}{60.0}{28.0}
\drawpath{60.0}{28.0}{56.0}{28.0}
\drawpath{56.0}{28.0}{46.0}{28.0}
\drawpath{60.0}{20.0}{56.0}{20.0}
\drawpath{56.0}{20.0}{46.0}{20.0}
\drawpath{66.0}{14.0}{84.0}{14.0}
\drawpath{60.0}{34.0}{62.0}{32.0}
\drawpath{64.0}{30.0}{66.0}{28.0}
\drawpath{60.0}{20.0}{62.0}{18.0}
\drawpath{64.0}{16.0}{66.0}{14.0}
\drawframebox{21.0}{38.0}{10.0}{12.0}{$R$}
\drawframebox{89.0}{10.0}{10.0}{12.0}{$S$}
\drawpath{46.0}{20.0}{38.0}{36.0}
\drawpath{46.0}{28.0}{42.0}{34.0}
\drawpath{42.0}{34.0}{38.0}{40.0}
\drawpath{38.0}{36.0}{36.0}{36.0}
\drawpath{36.0}{36.0}{34.0}{38.0}
\drawpath{34.0}{38.0}{36.0}{40.0}
\drawpath{36.0}{40.0}{38.0}{40.0}
\drawpath{46.0}{34.0}{44.0}{34.0}
\drawpath{38.0}{34.0}{26.0}{34.0}
\drawpath{46.0}{42.0}{28.0}{42.0}
\drawpath{28.0}{42.0}{26.0}{42.0}
\drawpath{4.0}{38.0}{4.0}{38.0}
\drawpath{4.0}{38.0}{4.0}{38.0}
\end{picture}
}
\end{center}

and hence to

\begin{center}
\centerline{\tt\setlength{\unitlength}{2.0pt}
\begin{picture}(98,50)
\thinlines
\drawpath{26.0}{44.0}{68.0}{44.0}
\drawpath{60.0}{34.0}{46.0}{34.0}
\drawpath{66.0}{28.0}{74.0}{28.0}
\drawpath{74.0}{28.0}{76.0}{24.0}
\drawpath{76.0}{24.0}{74.0}{20.0}
\drawpath{74.0}{20.0}{66.0}{20.0}
\drawpath{66.0}{20.0}{60.0}{14.0}
\drawpath{60.0}{14.0}{56.0}{14.0}
\drawpath{56.0}{14.0}{54.0}{10.0}
\drawpath{54.0}{10.0}{56.0}{6.0}
\drawpath{56.0}{6.0}{84.0}{6.0}
\drawpath{68.0}{44.0}{70.0}{42.0}
\drawpath{70.0}{42.0}{68.0}{40.0}
\drawpath{68.0}{40.0}{38.0}{40.0}
\drawframebox{21.0}{39.0}{10.0}{14.0}{$R$}
\drawpath{40.0}{34.0}{44.0}{34.0}
\drawpath{62.0}{32.0}{64.0}{30.0}
\drawpath{60.0}{20.0}{56.0}{20.0}
\drawpath{56.0}{20.0}{46.0}{20.0}
\drawpath{66.0}{14.0}{84.0}{14.0}
\drawpath{60.0}{34.0}{62.0}{32.0}
\drawpath{64.0}{30.0}{66.0}{28.0}
\drawpath{60.0}{20.0}{62.0}{18.0}
\drawpath{64.0}{16.0}{66.0}{14.0}
\drawframebox{89.0}{10.0}{10.0}{12.0}{$S$}
\drawpath{46.0}{20.0}{38.0}{36.0}
\drawpath{38.0}{36.0}{36.0}{36.0}
\drawpath{36.0}{36.0}{34.0}{38.0}
\drawpath{34.0}{38.0}{36.0}{40.0}
\drawpath{36.0}{40.0}{38.0}{40.0}
\drawpath{46.0}{34.0}{44.0}{34.0}
\drawpath{38.0}{34.0}{26.0}{34.0}
\drawpath{4.0}{38.0}{4.0}{38.0}
\drawpath{4.0}{38.0}{4.0}{38.0}
\end{picture}
}
\end{center}

This simplifies by duality to

\begin{center}
\centerline{\tt\setlength{\unitlength}{2.0pt}
\begin{picture}(98,40)
\thinlines
\drawpath{26.0}{34.0}{34.0}{34.0}
\drawpath{60.0}{34.0}{46.0}{34.0}
\drawpath{66.0}{28.0}{74.0}{28.0}
\drawpath{74.0}{28.0}{76.0}{24.0}
\drawpath{76.0}{24.0}{74.0}{20.0}
\drawpath{74.0}{20.0}{66.0}{20.0}
\drawpath{66.0}{20.0}{60.0}{14.0}
\drawpath{60.0}{14.0}{56.0}{14.0}
\drawpath{56.0}{14.0}{54.0}{10.0}
\drawpath{54.0}{10.0}{56.0}{6.0}
\drawpath{56.0}{6.0}{84.0}{6.0}
\drawpath{34.0}{34.0}{38.0}{26.0}
\drawpath{26.0}{24.0}{34.0}{24.0}
\drawpath{34.0}{24.0}{36.0}{26.0}
\drawframebox{21.0}{29.0}{10.0}{14.0}{$R$}
\drawpath{42.0}{34.0}{46.0}{34.0}
\drawpath{60.0}{34.0}{62.0}{32.0}
\drawpath{50.0}{20.0}{46.0}{20.0}
\drawpath{60.0}{20.0}{50.0}{20.0}
\drawpath{66.0}{14.0}{84.0}{14.0}
\drawpath{62.0}{32.0}{64.0}{30.0}
\drawpath{64.0}{30.0}{66.0}{28.0}
\drawpath{60.0}{20.0}{62.0}{18.0}
\drawpath{64.0}{16.0}{66.0}{14.0}
\drawpath{38.0}{28.0}{42.0}{34.0}
\drawframebox{89.0}{10.0}{10.0}{12.0}{$S$}
\drawpath{38.0}{26.0}{42.0}{20.0}
\drawpath{46.0}{34.0}{44.0}{34.0}
\drawpath{42.0}{20.0}{46.0}{20.0}
\drawpath{4.0}{28.0}{4.0}{28.0}
\drawpath{4.0}{28.0}{4.0}{28.0}
\end{picture}
}
\end{center}

It is now clear that repeating the argument using naturality and duality we obtain the result.
\qed

\section{A braided category of relations}\label{sec-trel}

\subsection{The definition of $\tr_G$}\label{subsec-defn-trel}
We will describe a braided modification of the category $\rel$ with a commutative Frobenius object.
\begin{definition}\label{def-trel}
Let $G$ be a group. The objects of $\tr_G$ are the formal powers of $G$, and
the arrows from $G^m$ to $G^n$
are relations $R$ from the set $G^m$ to the set $G^n$ satisfying:
\begin{itemize}
\item[1)] if $(x_1,...,x_m)R(y_1,...y_n)$ then also for all $g$ in $G$ \\
$(g^{-1}x_1g, ...,g^{-1}x_mg)R(g^{-1}y_1g, ... ,g^{-1}y_mg)$,
\item[2)] if $(x_1,...,x_m)R(y_1,...y_n)$ then $x_1...x_m(y_1...y_n)^{-1}\in Z(G)$ (the center of $G$).
\end{itemize}
Composition and identities are defined to be composition and identity of relations.
\end{definition}
It is straightforward to verify that $\tr_G$ is a category.
We introduce some useful notation.
Write $x = (x_1,...,x_m)$, $y = (y_1,...,y_n)$, and so on.
Write $\ol{x} = x_1x_2...x_m$ and for $g,h$ in $G$, as $g^h = hgh^{-1}$.
For $g$ in $G$ write $x^g = (x_1^g,x_2^g,...,x_m^g)$. Thus, $(\ol{x})^g = \ol{x^g}$, and
of course for any $x,y$ in $G^m\times G^n$,
$x^gy^g = (xy)^g$ where we write $xy$ for $(x_1,...,x_m,y_1,...,y_n)$.

\begin{theorem}\label{thm-trel-bmc}
$\tr_G$ is a braided strict monoidal category with tensor defined on objects by
$G^m\tn G^n = G^{m+n}$ and on arrows by product of relations. The twist 
$$\tau_{m,n}:G^m\tn G^n\to G^n\tn G^m$$
is the functional relation 
$$(x,y)\sim (y^{\ol{x}},x) $$
\end{theorem}

\pf
As noted above it is easy to show that identities and composites of relations satisfying
1) and 2) also satisfy 1) and 2), so $\tr_G$ is a category.  The monoidal structure of $\rel$ 
also restricts to $\tr_G$ since if $R:G^m\to G^t$ and $S:G^n\to G^u$ satisfy 1) and 2) then so also 
does $R\times S$. To see that $R\times S$ satisfies 1) notice that if $xRy$ and$zSw$ then for and $g\in G$, $x^gRy^g$ and $z^gSw^g$ and hence $(xz)^g(R\times S)(yw)^g$. To see that $R\times S$ satisfies 2) notice that, if $\ol{x}(\ol{y})^{-1}\in Z(G)$ and $\ol{z}(\ol{w})^{-1}\in Z(G)$, then
$\ol{xz}(\ol{yw})^{-1} = (\ol{x}) (\ol{z})((\ol{y}) (\ol{w}))^{-1} = (\ol{x}) (\ol{z})(\ol{w})^{-1}(\ol{y})^{-1}$.
But $\ol{z}(\ol{w})^{-1}\in Z(G)$, so $(\ol{x}) (\ol{z}) (\ol{w})^{-1}(\ol{y})^{-1} = %
(\ol{x})(\ol{y})^{-1}(\ol{z})(\ol{w})^{-1}$
and the latter is in $Z(G)$. 

We show that $B1$ holds for $\tw$ as defined. $B2$ is similar.

First note that $\tw_{m,n+p}(xyz)=(yz)^{\ol{x}}x$.
Further $(\tw_{m,n}\tn 1_{G^p})(xyz)= y^{\ol{x}}xz$
while $(1_{G^n}\tn \tw_{n,p})(xyz)= xz^{\ol{y}}y \}$.
Thus $$(1_{G^n}\tn \tw_{n,p})(\tw_{m,n}\tn 1_{G^p})(xyz) = (1_{G^n}\tn \tw_{n,p})((y^{\ol{x}})xz) = (y^{\ol{x}})(z^{\ol{x}})x=\tw_{m,n+p}(xyz).$$ 

Lastly we need to show that $\tw_{m,n}:G^m\times G^n\to G^n\times G^m$ is natural. This amounts to two conditions. Consider $R:G^p\to G^m$ and $S:G^q\to G^n$ in $\tr_G$. 
The first condition for naturality is that 
$$\tw_{m,n} (R\tn 1_{G^n})= (1_{G^n}\tn R)\tw_{p,n}:G^{p+n}\to G^{n+m}.$$ 
But $xyzw$ ($x\in G^p$, $y\in  G^n$, $z\in G^n$, $w\in G^m $) belongs to the left-hand side iff $x R w$ and $z=y^{\ol{w}}$, whereas $xyzw$ belongs to the right-hand side iff $xRw$ and $z=y^{\ol{x}} $. But condition 2) implies that if $xRw$ then for any $y$ it follows that $y^{\ol{x}}=y^{\ol{w}}$, and hence the result. 

The second condition for naturality is that
$$(\tw_{m,n})(1_{G^m}\tn S)=(S\tn 1_{G^m})(\tw_{m,q}):G^{m+q}\to G^{n+m}. $$ 
But $xyzw$ ($x\in G^m$, $y\in  G^q$, $z\in G^n$, $w\in G^m $) belongs to the left-hand side iff
$x=w$ and $yS(z^{\ol{x}^{-1}}) $, whereas $xyzw$ belongs to the right-hand side iff $x=w$ and $y^{\ol{x}}Sz$. Condition 1) implies the result.
\qed

\begin{rem}
Notice that a relation in $\tr_G$ from $I$ to $G\times G$ is just a subset of $G\times G$ closed under conjugation by elements of $G$ and whose elements $(x,y)$ satisfy $xy\in centre(G)$. Further a relation from $I$ to $I$ is either the empty set or the one-point set. 

Notice also that if the group $G$ is abelian the conditions $(1)$ and $(2)$ of the definition \ref{def-trel} are trivially true.
\end{rem}

\subsubsection{The commutative Frobenius structure on $G$}\label{subsubsec-frob-alg-trel}

The commutative Frobenius structure on the object $G$ of $\tr_G$ mentioned above is as follows:
$\nabla$ is a function, namely the multiplication of the group $G$, $n: I \to G$ is also a function,
the identity of the group; $\Delta$ is the opposite relation of $\nabla$, $e$ is the opposite relation
of $n$. 

Notice that $\eta$ is the relation $*\sim (x,x^{-1})$, and $\eps$ is the opposite relation of $\eta$.

It is straightforward to check that these relations belong to $\tr_G$. We will just check one of the Frobenius equations, namely that $$(1_G\tn \nabla)(\Delta\tn 1_G)=\Delta\nabla :G\times G\to G\times G.$$ 
If $g$, $h$, $p$, $q$ are in $G$ then $(g,h,p,q)$ belongs to the left-hand relation if 
there is a $r\in G$ such that $g=pr$ and $rh=q$. But this is the same as saying that $p^{-1}g=qh^{-1}$ or $gh=pq$ which is exactly the condition for $(g,h,p,q)$ to be in the right-hand relation. 

\subsection{Proving circuits distinct in $\tr$}\label{subsec-proving-distinct}
In this section we discuss the possibility of distinguishing various
tangled circuits, including the analogue of knots, closed circuits, that is, circuits 
from the one-point set $I$ to $I$, by looking in $\tr_G$.

\subsubsection{Example}\label{subsubsec-four-wires}
First an example where two circuits may be distinguished in $\tr_{S_3}$, where $S_3$ is the symmetric group on three letters. The circuits are:

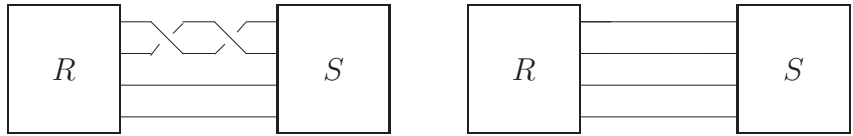
\begin{center}
{\tt\setlength{\unitlength}{3.0pt}
\begin{picture}(114,24)
\thinlines
\drawframebox{11.0}{12.0}{14.0}{16.0}{$R$}
\drawpath{18.0}{18.0}{22.0}{18.0}
\drawpath{18.0}{14.0}{22.0}{14.0}
\drawpath{22.0}{18.0}{26.0}{14.0}
\drawpath{26.0}{14.0}{30.0}{14.0}
\drawpath{26.0}{18.0}{30.0}{18.0}
\drawpath{30.0}{18.0}{34.0}{14.0}
\drawpath{34.0}{14.0}{38.0}{14.0}
\drawpath{34.0}{18.0}{38.0}{18.0}
\drawpath{18.0}{10.0}{38.0}{10.0}
\drawpath{18.0}{6.0}{38.0}{6.0}
\drawframebox{45.0}{12.0}{14.0}{16.0}{$S$}
\drawpath{21.97}{13.75}{23.04}{15.15}
\drawpath{24.79}{16.54}{26.02}{17.77}
\drawpath{30.05}{13.93}{31.29}{15.15}
\drawpath{32.88}{16.72}{33.93}{17.93}
\drawframebox{69.0}{12.0}{14.0}{16.0}{$R$}
\drawpath{76.0}{18.0}{80.0}{18.0}
\drawpath{76.0}{14.0}{80.0}{14.0}
\drawpath{76.0}{18.0}{96.0}{18.0}
\drawpath{84.0}{14.0}{88.0}{14.0}
\drawpath{76.0}{14.0}{96.0}{14.0}
\drawpath{92.0}{14.0}{96.0}{14.0}
\drawpath{76.0}{10.0}{96.0}{10.0}
\drawpath{76.0}{6.0}{96.0}{6.0}
\drawframebox{103.0}{12.0}{14.0}{16.0}{$S$}
\end{picture}
}
\end{center}

\pf
Let each of $R$ and $S$ be the set of conjugates of $u=(12,13,23,13)$ under the action of $G$ (not $G\times G\times G\times G$).
Notice that $(12)(13)(23)(13)$ is the identity.
The second circuit evaluates as the one point set.

The first circuit evaluates instead as the empty set since the braid in the first circuit relates $(12,13,23,13)$ in $R$ to $(13,23,23,13)$ which is not in the conjugacy class of $u$ since the second and third elements are equated by the braid. 
 
\qed
Notice that a similar argument using the symmetric group $S_{3}$ works for two components joined by $n>3$ wires, the first two of which are tangled.

\subsubsection{Example}\label{subsubsec-intro-circs}
We will see that the first two circuits in section \ref{subsec-tangled-circuits} can also be shown distinct in $\tr_{S_3}$.  It is clearly sufficient to show the following circuits distinct:

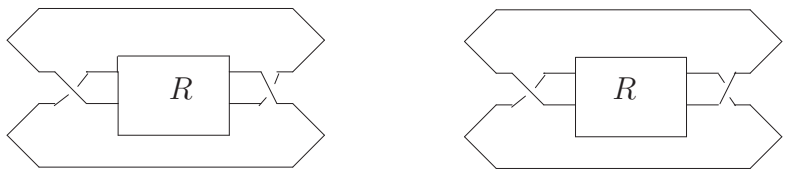
\begin{center}
{\tt\setlength{\unitlength}{3.0pt}
\begin{picture}(106,28)
\thinlines
\drawpath{4.13}{20.02}{8.13}{24.02}
\drawpath{8.13}{24.02}{40.13}{24.02}
\drawpath{40.13}{24.02}{44.13}{20.02}
\drawpath{4.13}{20.02}{8.13}{16.02}
\drawpath{8.13}{16.02}{10.13}{16.02}
\drawpath{10.13}{16.02}{14.12}{12.02}
\drawpath{14.12}{12.02}{18.12}{12.02}
\drawpath{18.12}{12.02}{18.12}{16.02}
\drawpath{18.03}{15.97}{18.03}{17.95}
\drawpath{32.13}{16.02}{32.13}{8.02}
\drawpath{18.12}{18.02}{32.13}{18.02}
\drawpath{18.12}{8.02}{18.12}{12.02}
\drawpath{32.13}{18.02}{32.13}{16.02}
\drawpath{18.12}{16.02}{14.12}{16.02}
\drawpath{14.12}{16.02}{14.12}{16.02}
\drawpath{14.29}{15.77}{13.04}{14.18}
\drawpath{10.06}{11.76}{11.84}{13.14}
\drawpath{10.13}{12.02}{8.13}{12.02}
\drawpath{8.13}{12.02}{4.13}{8.02}
\drawpath{4.13}{8.02}{8.13}{4.01}
\drawpath{8.13}{4.01}{40.13}{4.01}
\drawpath{40.13}{4.01}{44.13}{8.02}
\drawpath{44.13}{8.02}{40.13}{12.02}
\drawpath{44.13}{20.02}{40.13}{16.02}
\drawpath{40.13}{16.02}{38.13}{16.02}
\drawpath{40.13}{12.02}{38.13}{12.02}
\drawpath{36.13}{16.02}{38.13}{12.02}
\drawpath{36.13}{16.02}{32.13}{16.02}
\drawpath{32.13}{12.02}{36.13}{12.02}
\drawpath{18.12}{8.02}{32.13}{8.02}
\drawpath{36.79}{12.77}{35.9}{11.76}
\drawpath{38.38}{15.77}{37.86}{14.52}
\drawpath{61.79}{19.84}{65.79}{23.84}
\drawpath{65.79}{23.84}{97.79}{23.84}
\drawpath{97.79}{23.84}{101.79}{19.84}
\drawpath{61.79}{19.84}{65.79}{15.84}
\drawpath{65.79}{15.84}{67.79}{15.84}
\drawpath{67.79}{15.84}{71.79}{11.84}
\drawpath{71.79}{11.84}{75.79}{11.84}
\drawpath{75.79}{11.84}{75.79}{15.84}
\drawpath{75.79}{15.84}{75.79}{17.84}
\drawpath{89.79}{15.84}{89.79}{7.84}
\drawpath{75.79}{17.84}{89.79}{17.84}
\drawpath{75.79}{7.84}{75.79}{11.84}
\drawpath{89.79}{17.84}{89.79}{15.84}
\drawpath{75.79}{15.84}{71.79}{15.84}
\drawpath{71.79}{15.84}{71.79}{15.84}
\drawpath{71.94}{15.59}{70.72}{14.02}
\drawpath{67.73}{11.56}{69.48}{12.97}
\drawpath{67.79}{11.84}{65.79}{11.84}
\drawpath{65.79}{11.84}{61.79}{7.84}
\drawpath{61.79}{7.84}{65.79}{3.82}
\drawpath{65.79}{3.82}{97.79}{3.82}
\drawpath{97.79}{3.82}{101.79}{7.84}
\drawpath{101.79}{7.84}{97.79}{11.84}
\drawpath{101.79}{19.84}{97.79}{15.84}
\drawpath{97.79}{15.84}{95.79}{15.84}
\drawpath{97.79}{11.84}{95.79}{11.84}
\drawpath{96.0}{16.0}{94.0}{12.0}
\drawpath{93.79}{15.84}{89.79}{15.84}
\drawpath{89.79}{11.84}{93.79}{11.84}
\drawpath{75.79}{7.84}{89.79}{7.84}
\drawpath{93.83}{15.75}{94.54}{14.68}
\drawpath{95.25}{12.77}{95.94}{11.75}
\drawcenteredtext{26.0}{14.0}{$R$}
\drawcenteredtext{82.0}{14.0}{$R$}
\end{picture}
}
\end{center}

Take $R$ to be the following subset of $(S_3)^2\times (S_3)^2$: the conjugacy class of the element 
$((12,13),(12,13))$. Then the first circuit evaluates as $\emptyset$ and the second as the one-point set. 
\subsubsection{Example}\label{subsubsec-three-wires-equal-in-TRel}
Next an example of two circuits which we believe are distinct in $\tcd_M$ but
are always equal in $\tr_G$. For any group $G$, $\tr$ cannot distinguish them.

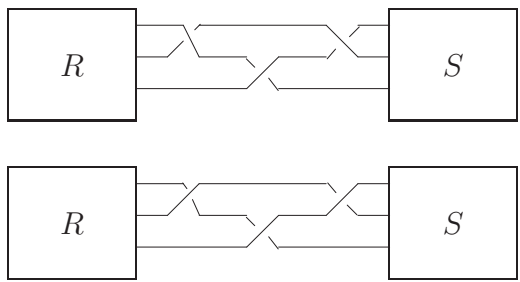
\begin{center}

{\tt\setlength{\unitlength}{3.0pt}
\begin{picture}(72,42)
\thinlines
\drawframebox{12.0}{31.0}{16.0}{14.0}{$R$}
\drawpath{20.0}{36.0}{26.0}{36.0}
\drawpath{26.0}{36.0}{28.0}{32.0}
\drawpath{28.0}{32.0}{34.0}{32.0}
\drawpath{20.0}{32.0}{24.0}{32.0}
\drawpath{24.0}{32.0}{26.0}{34.0}
\drawpath{28.0}{36.0}{44.0}{36.0}
\drawpath{28.13}{36.12}{27.43}{35.41}
\drawpath{20.0}{28.0}{34.0}{28.0}
\drawpath{34.0}{28.0}{38.0}{32.0}
\drawpath{38.0}{32.0}{44.0}{32.0}
\drawpath{38.0}{28.0}{52.0}{28.0}
\drawpath{33.93}{31.75}{35.0}{30.7}
\drawpath{36.75}{29.31}{37.97}{27.75}
\drawpath{44.0}{36.0}{48.0}{32.0}
\drawpath{48.0}{32.0}{52.0}{32.0}
\drawpath{48.0}{36.0}{52.0}{36.0}
\drawpath{45.18}{32.97}{44.13}{31.75}
\drawpath{48.18}{35.77}{46.95}{34.54}
\drawframebox{60.0}{31.0}{16.0}{14.0}{$S$}
\drawframebox{12.0}{11.0}{16.0}{14.0}{$R$}
\drawpath{20.0}{16.0}{26.0}{16.0}
\drawpath{24.0}{12.0}{28.0}{16.0}
\drawpath{28.0}{12.0}{34.0}{12.0}
\drawpath{20.0}{12.0}{24.0}{12.0}
\drawpath{44.0}{12.0}{48.0}{16.0}
\drawpath{28.0}{16.0}{44.0}{16.0}
\drawpath{27.15}{13.99}{28.02}{12.06}
\drawpath{20.0}{8.0}{34.0}{8.0}
\drawpath{34.0}{8.0}{38.0}{12.0}
\drawpath{38.0}{12.0}{44.0}{12.0}
\drawpath{38.0}{8.0}{52.0}{8.0}
\drawpath{33.93}{11.75}{35.0}{10.7}
\drawpath{36.75}{9.31}{37.97}{7.75}
\drawpath{44.2}{16.08}{45.25}{14.86}
\drawpath{48.0}{12.0}{52.0}{12.0}
\drawpath{48.0}{16.0}{52.0}{16.0}
\drawpath{46.65}{13.29}{47.88}{12.06}
\drawpath{25.91}{16.08}{26.45}{15.38}
\drawframebox{60.0}{11.0}{16.0}{14.0}{$S$}
\end{picture}
}
\end{center}
\pf
Suppose $(x,y,z)$ is an element of component $R$.  Notice that since $xyz$ is in the centre $xyz=yzx=zxy$.  The braid between the two components in the first circuit relates $(x,y,z)$ to $$u=(xyx^{-1}zxy^{-1}x^{-1},xyx^{-1},z^{-1}xz)=(xyx^{-1}zxy^{-1}x^{-1},z^{-1}yz,z^{-1}xz)$$ since $yzx=zxy$. Instead the braid in the second circuit relates $(x,y,z)$ to 
$v=(z,z^{-1}yz,z^{-1}y^{-1}xyz)=(z,xyx^{-1},x)$ since $z^{-1}y^{-1}xyz=z^{-1}y^{-1}yzx=x$ and $zxy=yzx$. But $xzuz^{-1}x^{-1}=v$ since $$xzxyx^{-1}zxy^{-1}x^{-1}z^{-1}x^{-1}=xyzxx^{-1}zxx^{-1}z^{-1}y^{-1}x^{-1}=xyzy^{-1}x^{-1}=z$$ and hence $u$ and $v$ are conjugate. 
Since $S$ is closed under conjugacy, the element $(x,y,z)$ gives rise to an element of the first circuit if and only if it does for the second circuit. Since this is true  for any $(x,y,z)$ the two circuits are equal in $\tr_G$.
\qed
\subsubsection{Example}\label{subsubsec-three-wires-any-braid-equal-in-TRel}
In fact the last example is general for three wires. \emph{The circuit obtained by composing in $\tr_G$ any two  two components $R:I\to G^3$ and $S:G^3\to I$ with a braiding in between depends only on the permutation, not the braiding}.

\pf
Suppose $(x,y,z)\in R$ then $xyz\in centre(G)$ and hence $xyx^{-1}=z^{-1}yz$, $yzy^{-1}=x^{-1}zx$ and $zxz^{-1}=y^{-1}xy$. Consider two composites  $R$ composed with $\tau \tn 1$ and $R$ composed with $\tau^{-1} \tn 1$. Consider $(x,y,z)\in R$. We will show that these two composites associate $(x,y,z)$ with conjugate triples. Repeating this we see that the argument given in the above example can be applied, showing that in a composite $\tau$ and $\tau^{-1}$ are interchangeable.

In the first composite $(x,y,z)$ is related to $u=(xyx^{-1},x,z)=(z^{-1}yz,y,z)$. In the second composite $(x,y,z)$ is related to $(y,y^{-1}xy,z)$. It is immediate that $zuz^{-1}=v$.   
\qed 

Of course different permutations can be distinguished even in $Rel$.

\subsubsection{Example}\label{subsubsec-two-houses-four-utilities}
Another two circuits we can distinguish in $\tr_{S_3}$:

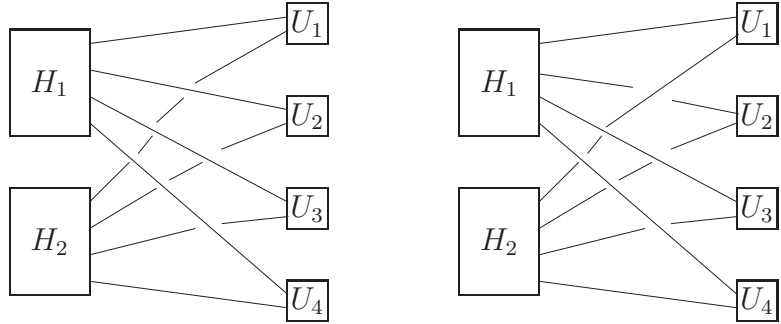
\begin{center}
{\tt\setlength{\unitlength}{2.5pt}
\begin{picture}(124,56)
\thinlines
\drawframebox{9.81}{40.63}{12.0}{16.0}{$H_1$}
\drawframebox{9.81}{16.63}{12.0}{16.0}{$H_2$}
\drawframebox{48.81}{49.63}{6.0}{6.0}{$U_1$}
\drawframebox{48.81}{35.63}{6.0}{6.0}{$U_2$}
\drawframebox{48.81}{21.63}{6.0}{6.0}{$U_3$}
\drawframebox{48.81}{7.63}{6.0}{6.0}{$U_4$}
\drawpath{15.81}{46.63}{45.81}{50.63}
\drawpath{15.81}{42.63}{45.81}{36.63}
\drawpath{15.81}{38.63}{45.81}{22.63}
\drawpath{15.81}{34.63}{45.81}{8.63}
\drawpath{15.81}{22.63}{21.81}{28.63}
\drawpath{15.81}{18.63}{25.81}{24.63}
\drawpath{25.81}{24.63}{25.81}{24.63}
\drawpath{15.81}{14.63}{31.81}{18.63}
\drawpath{15.81}{10.63}{45.81}{6.63}
\drawpath{45.81}{48.63}{31.81}{40.63}
\drawpath{45.81}{34.63}{35.81}{30.63}
\drawpath{30.06}{38.38}{26.02}{34.72}
\drawpath{25.31}{32.63}{23.2}{30.0}
\drawpath{32.88}{28.61}{27.95}{26.0}
\drawpath{45.88}{20.43}{35.86}{19.38}
\drawframebox{77.8}{40.63}{12.0}{16.0}{$H_1$}
\drawframebox{77.81}{16.63}{12.0}{16.0}{$H_2$}
\drawframebox{116.81}{49.63}{6.0}{6.0}{$U_1$}
\drawframebox{116.81}{35.63}{6.0}{6.0}{$U_2$}
\drawframebox{116.81}{21.63}{6.0}{6.0}{$U_3$}
\drawframebox{116.81}{7.63}{6.0}{6.0}{$U_4$}
\drawpath{83.81}{46.63}{113.81}{50.63}
\drawpath{94.0}{34.0}{114.0}{48.0}
\drawpath{83.81}{38.63}{113.81}{22.63}
\drawpath{83.81}{34.63}{113.81}{8.63}
\drawpath{83.81}{22.63}{89.81}{28.63}
\drawpath{83.81}{18.63}{93.81}{24.63}
\drawpath{93.81}{24.63}{93.81}{24.63}
\drawpath{83.81}{14.63}{99.81}{18.63}
\drawpath{83.81}{10.63}{113.81}{6.63}
\drawpath{84.0}{42.0}{98.0}{40.0}
\drawpath{113.81}{34.63}{103.81}{30.63}
\drawpath{104.0}{38.0}{114.0}{36.0}
\drawpath{93.31}{32.63}{91.2}{30.0}
\drawpath{100.88}{28.61}{95.95}{26.0}
\drawpath{113.88}{20.43}{103.86}{19.38}
\end{picture}
}
\end{center}

\pf
Replace each of the four components $U_1$,$U_2$,$U_3$,$U_4$ by $\eps$. Let $R$ be the conjugacy class of $(12,13,23,13)$. The wires of the first circuit relate this element to $u=(12,23,12,13)$, and of the second circuit to $v=(13,12,12,13)$.Clearly $u$ and $v$ are not conjugate, and hence we can choose $S$ so that the two circuits evaluate differently in $\tr_{S_3}$.
\qed
\subsubsection{Example}\label{subsubsec-forks}
The following two circuits can be distinguished in $\tr_{S_3}$.
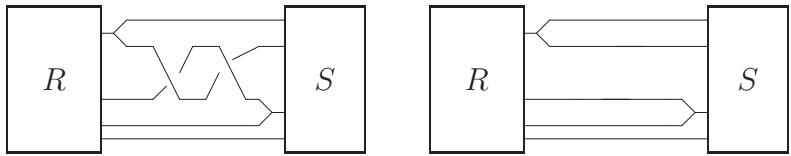
\begin{center}
{\tt\setlength{\unitlength}{2.5pt}
\begin{picture}(126,30)
\thinlines
\drawframebox{11.0}{15.0}{14.0}{22.0}{$R$}
\drawframebox{52.0}{15.0}{12.0}{22.0}{$S$}
\drawpath{18.0}{6.0}{46.0}{6.0}
\drawpath{18.0}{22.0}{20.0}{22.0}
\drawpath{20.0}{22.0}{22.0}{24.0}
\drawpath{22.0}{24.0}{46.0}{24.0}
\drawpath{20.0}{22.0}{22.0}{20.0}
\drawpath{22.0}{20.0}{26.0}{20.0}
\drawpath{46.0}{10.0}{44.0}{10.0}
\drawpath{44.0}{10.0}{42.0}{8.0}
\drawpath{42.0}{8.0}{40.0}{8.0}
\drawpath{44.0}{10.0}{42.0}{12.0}
\drawpath{42.0}{12.0}{40.0}{12.0}
\drawpath{40.0}{8.0}{18.0}{8.0}
\drawpath{26.0}{20.0}{30.0}{12.0}
\drawpath{30.0}{12.0}{34.0}{12.0}
\drawpath{18.0}{12.0}{26.0}{12.0}
\drawpath{32.0}{20.0}{36.0}{20.0}
\drawpath{36.0}{20.0}{40.0}{12.0}
\drawpath{46.0}{20.0}{42.0}{20.0}
\drawpath{42.0}{20.0}{38.0}{18.0}
\drawpath{36.0}{16.0}{34.0}{12.0}
\drawpath{32.0}{20.0}{30.0}{16.0}
\drawpath{28.0}{14.0}{26.0}{12.0}
\drawframebox{75.0}{15.0}{14.0}{22.0}{$R$}
\drawframebox{116.0}{15.0}{12.0}{22.0}{$S$}
\drawpath{82.0}{6.0}{110.0}{6.0}
\drawpath{82.0}{22.0}{84.0}{22.0}
\drawpath{84.0}{22.0}{86.0}{24.0}
\drawpath{86.0}{24.0}{110.0}{24.0}
\drawpath{84.0}{22.0}{86.0}{20.0}
\drawpath{86.0}{20.0}{90.0}{20.0}
\drawpath{110.0}{10.0}{108.0}{10.0}
\drawpath{108.0}{10.0}{106.0}{8.0}
\drawpath{106.0}{8.0}{104.0}{8.0}
\drawpath{108.0}{10.0}{106.0}{12.0}
\drawpath{106.0}{12.0}{104.0}{12.0}
\drawpath{104.0}{8.0}{82.0}{8.0}
\drawpath{90.0}{20.0}{106.0}{20.0}
\drawpath{94.0}{12.0}{98.0}{12.0}
\drawpath{82.0}{12.0}{90.0}{12.0}
\drawpath{96.0}{20.0}{100.0}{20.0}
\drawpath{90.0}{12.0}{104.0}{12.0}
\drawpath{110.0}{20.0}{106.0}{20.0}
\end{picture}
}
\end{center}
\pf
Take $R$ to be the conjugacy class of $(12,13,23,13)$ and $S$ the conjugacy class of 
$((),13,(),13)$. The first circuit evaluates as  the one-point set and the second as $\emptyset$.
\qed
\section{A braided category of spans}\label{sec-tspan}

In this section we begin to extend the previous sections with a modification of the category $\spa$ of spans of sets
with a braiding for some spans.

\begin{definition}\label{def-tspan}
Let $G$ be a group. The objects of $\tsp_G$ are the formal powers of $G$, and
an arrow from $G^m$ to $G^n$
is an isomorphism class of spans in sets, $G^m \toleft^{\delta_0} S \to^{\delta_1} G^n$, from the set $G^m$ to the set $G^n$
such that there exist a function $G\times S \to S$ of $G$ written $(g,s)\mapsto gs$ yielding a bijection for each $g\in G$, and
satisfying:
\begin{itemize}
\item[1)] if $\delta_0(s) = (x_1,...,x_m)$ and $\delta_1(s) = (y_1,...,y_n)$ then %
$\delta_0(gs) = (x_1^g,...,x_m^g)$  and $\delta_1(gs) = (y_1^g, ... ,y_m^g)$ for all $g$ in $G$,
\item[2)] if $\delta_0(s) = (x_1,...,x_m)$ and $\delta_1(s) = (y_1,...,y_n)$  then $x_1...x_m(y_1...y_n)^{-1}\in Z(G)$.
\end{itemize}
Composition and identities are composition and identity of spans.
\end{definition}

It is straightforward that $\tsp_G$ is a category.
Like $\tr_G$ it has the structure of a braided
strict monoidal category.

\begin{theorem}\label{thm-tspan-bmc}
$\tsp_G$ is braided strict monoidal with tensor defined by $G^m\tn G^n = G^{m+n}$ and twist
$\tw_{m,n} : G^m\tn G^n \to G^n\tn G^m$ is the span
determined by the function $\delta_1$ where:
$$ \delta_1(x_1,...,x_m,y_1,...,y_n) = (y_1^{\ol{x}},...,y_n^{\ol{x}},x_1,...,x_m) $$
where $\ol{x} = x_1x_2...x_m$ and $y^x =xyx^{-1}$. 
\end{theorem}

\pf
This is similar to Theorem \ref{thm-trel-bmc}. We use the same notation as above.
As noted, it is easy to show that identities and composites of spans satisfying conditions 1) and 2)
also satisfy 1) and 2), so $\tsp_G$ is a category.

To see that $\tn$ is a functor recall that product of spans defines a tensor functor on the
category $\spa$ of spans. It remains to show that $\tsp_G$ is closed under $\tn$. Suppose $R : G^m \to G^t$
and $S : G^n \to G^u$.
If $x = \delta_0(r), y = \delta_1(r)$ and $z = \delta_0(s), w \delta_1(s)$,
then for any $g$, $x^g = \delta_0(gr), y^g = \delta_1(gr)$ and $z^g = \delta_0(s), w^g = \delta_1(s)$,
whence $(xz)^g =\delta_0(gr,gs), (yw)^g = \delta_1(gr,gs)$, so taking $g(r,s)$ to be $(gr,gs)$ condition 1) is satisfied.
For $x, y, z, w$ as defined, condition 2) follows exactly as in Theorem \ref{thm-trel-bmc}.

The associative and unitary properties for $\tn$ in $\tsp_G$ are immediate from the same properties in $\spa$.

We show that $B1$ holds for $\tw$ as defined. $B2$ is similar.

Since the twists and identities are defined by functions, the span composition is obtained by composing
functions and we calculate:
\begin{align}
(1_n\tn\tw_{m,p})(\tw_{m,n}\tn 1_p)(xyz) & = (1_n\tn\tw_{m,p})(y^{\ol{x}}xz) \\
& = y^{\ol{x}}z^{\ol{x}}x\\
& = (yz)^{\ol{x}}x\\
& = \tw_{m,n+p}(xyz)\\
\end{align}
which proves $B1$.

As in the case of $\tr_G$ the conditions 1) and 2) assure the naturality of $\tw$.
\qed

\subsubsection{A commutative Frobenius structure on $G$}\label{subsubsec-frob-alg-tspan}
As for $\tr_G$ and using the same functions viewed as spans, $G$ has the structure of a commutative
Frobenius algebra in $\tsp_G$. Consequently:

\begin{cor}
There is a unique braided strict monoidal functor $$\tang\to \tsp_G$$ taking the generating object to $G$ and
the structure maps of
$\tang$ to the corresponding arrows in $\tsp_G$.
\end{cor}
\subsection{Knot colourings}\label{subsec-knot-colourings}
The description of $\tsp_G$ makes it clear that there is a faithful monoidal functor
$$\tsp_G\to \spa(\set).$$ The following composite of monoidal functors we have described we denote as $colourings$:
$$colourings_G :\tang \to \tcd \to \tsp_G \to \spa (\set ).$$ 
$colourings_G$ takes the generating object $X$ of $\tang$ to the underlying set of $G$, and takes $\eps_X$ to the span $G\times G\leftarrow \{(x,y):xy=1\}\to I$, $\eta_X$ to $I\leftarrow \{(x,y):xy=1\}\to G\times G$ and $\tau_X$ to $(x,y)\leftarrow (x,y)\mapsto (xyx^{-1},x)$.

\begin{thm} (J. Armstrong \cite{accn})
 If $K$ is a knot then $colourings_{G}$ is the set of colourings of $K$ in the group $G$.
\end{thm}

\begin{rem}
Because of the faithfulness of the functor $\tsp_G\to \spa(\set)$  the calculation of the set of colourings of a knot may be done equally in $\tsp_G$ or $\spa(\set)$. The advantage of introducing $\tsp_G$ as we do is that $\tsp_G$ has the same structure as $\tang$ (braided monoidal with a tangle algebra) whereas $\spa (\set)$ does not.
\end{rem}

\subsubsection{Colourings of a trefoil}\label{subsubsec-colour-trefoil}
We will calculate the colourings of a trefoil in the dihedral group $D_3$ to allow us to introduce notation and indicate relations with other work. One expression for a trefoil in $\tang$ is 
$$(\eps\tn\eps)(1\tn \tau\tn 1)(1\tn 1\tn \tau^{-1}(1\tn \tau\tn 1)(\eta\tn\eta).$$

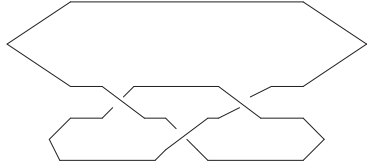
\begin{center}
{\tt\setlength{\unitlength}{2.0pt} 
\begin{picture}(76,38)
\thinlines
\drawpath{16.0}{34.0}{4.0}{26.0}
\drawpath{16.0}{34.0}{60.0}{34.0}
\drawpath{60.0}{34.0}{72.0}{26.0}
\drawpath{72.0}{26.0}{60.0}{18.0}
\drawpath{60.0}{18.0}{54.0}{18.0}
\drawpath{4.0}{26.0}{16.0}{18.0}
\drawpath{16.0}{18.0}{22.0}{18.0}
\drawpath{22.0}{18.0}{30.0}{12.0}
\drawpath{30.0}{12.0}{34.0}{12.0}
\drawpath{22.0}{12.0}{24.0}{14.0}
\drawpath{26.0}{16.0}{28.0}{18.0}
\drawpath{28.0}{18.0}{34.0}{18.0}
\drawpath{22.0}{12.0}{16.0}{12.0}
\drawpath{42.0}{12.0}{34.0}{6.0}
\drawpath{34.0}{12.0}{36.0}{10.0}
\drawpath{38.0}{8.0}{40.0}{6.0}
\drawpath{40.0}{6.0}{42.0}{4.0}
\drawpath{42.0}{12.0}{44.0}{12.0}
\drawpath{34.0}{18.0}{44.0}{18.0}
\drawpath{44.0}{18.0}{52.0}{12.0}
\drawpath{44.0}{12.0}{48.0}{14.0}
\drawpath{50.0}{16.0}{54.0}{18.0}
\drawpath{52.0}{12.0}{60.0}{12.0}
\drawpath{34.0}{6.0}{32.0}{4.0}
\drawpath{32.0}{4.0}{14.0}{4.0}
\drawpath{14.0}{4.0}{12.0}{8.0}
\drawpath{12.0}{8.0}{16.0}{12.0}
\drawpath{42.0}{4.0}{60.0}{4.0}
\drawpath{60.0}{4.0}{64.0}{8.0}
\drawpath{64.0}{8.0}{60.0}{12.0}
\end{picture}
}
\end{center}

It is convenient to represent the arrows in this expression as components as follows:

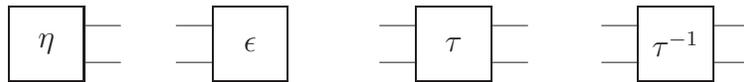
\begin{center}
{\tt\setlength{\unitlength}{3.5pt}
\begin{picture}(88,16)
\thinlines
\drawframebox{8.0}{8.0}{8.0}{8.0}{$\eta$}
\drawpath{12.0}{10.0}{16.0}{10.0}
\drawframebox{30.0}{8.0}{8.0}{8.0}{$\epsilon$}
\drawpath{12.0}{6.0}{16.0}{6.0}
\drawpath{22.0}{10.0}{26.0}{10.0}
\drawpath{22.0}{6.0}{26.0}{6.0}
\drawframebox{52.0}{8.0}{8.0}{8.0}{$\tau$}
\drawpath{44.0}{10.0}{48.0}{10.0}
\drawpath{44.0}{6.0}{48.0}{6.0}
\drawpath{56.0}{10.0}{60.0}{10.0}
\drawpath{56.0}{6.0}{60.0}{6.0}
\drawframebox{76.0}{8.0}{8.0}{8.0}{$\tau^{-1}$}
\drawpath{68.0}{10.0}{72.0}{10.0}
\drawpath{68.0}{6.0}{72.0}{6.0}
\drawpath{80.0}{10.0}{84.0}{10.0}
\drawpath{80.0}{6.0}{84.0}{6.0}
\end{picture}
}

\end{center}

Then the trefoil may be written as the circuit diagram:

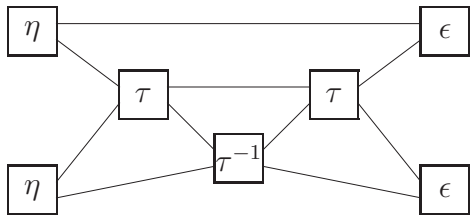
\begin{center}
{\tt\setlength{\unitlength}{3.0pt}
\begin{picture}(66,34)
\thinlines
\drawframebox{7.0}{27.0}{6.0}{6.0}{$\eta$}
\drawframebox{7.0}{7.0}{6.0}{6.0}{$\eta$}
\drawframebox{21.0}{19.0}{6.0}{6.0}{$\tau$}
\drawframebox{33.0}{11.0}{6.0}{6.0}{$\tau^{-1}$}
\drawframebox{45.0}{19.0}{6.0}{6.0}{$\tau$}
\drawframebox{59.0}{27.0}{6.0}{6.0}{$\eps$}
\drawframebox{59.0}{7.0}{6.0}{6.0}{$\eps$}
\drawpath{10.0}{28.0}{56.0}{28.0}
\drawpath{10.0}{26.0}{18.0}{20.0}
\drawpath{18.0}{18.0}{10.0}{8.0}
\drawpath{10.0}{6.0}{30.0}{10.0}
\drawpath{24.0}{18.0}{30.0}{12.0}
\drawpath{24.0}{20.0}{42.0}{20.0}
\drawpath{36.0}{12.0}{42.0}{18.0}
\drawpath{56.0}{26.0}{48.0}{20.0}
\drawpath{48.0}{18.0}{56.0}{8.0}
\drawpath{36.0}{10.0}{56.0}{6.0}
\end{picture}
}
\end{center}

The evaluation of the expression for the trefoil in $\spa(\set)$ is a limit of the diagram in $\set$ formed by taking for each wire in the diagram the set $G$ and for each component the pair of arrows constituting its span of sets (see \cite{RSW08} for the relation between limits in $\bC$ and expressions in $\spa (\bC)$). An element of this limit is a tuple of elements of $G$ one for each wire, satisfying the conditions of the components. Each of the components $\eta$, $\eps$, $\tau$, $\tau^{-1}$ is actually a relation from its domain to codomain, that is a subset of products of groups given by equational conditions. 

It is convenient to refine the pictures of the component to include the conditions as follows:
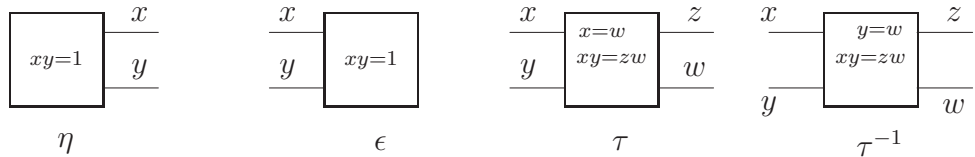
\begin{center}
{\tt\setlength{\unitlength}{3.5pt}
\begin{picture}(112,26)
\thinlines
\drawframebox{9.0}{15.0}{10.0}{10.0}{$\scriptstyle{xy=1}$}
\drawpath{14.0}{18.0}{20.0}{18.0}
\drawpath{14.0}{12.0}{20.0}{12.0}
\drawcenteredtext{18.0}{20.0}{$x$}
\drawcenteredtext{18.0}{14.0}{$y$}
\drawcenteredtext{10.0}{6.0}{$\eta$}
\drawframebox{43.0}{15.0}{10.0}{10.0}{$\scriptstyle{xy=1}$}
\drawpath{32.0}{18.0}{38.0}{18.0}
\drawpath{32.0}{12.0}{38.0}{12.0}
\drawcenteredtext{34.0}{20.0}{$x$}
\drawcenteredtext{34.0}{14.0}{$y$}
\drawframebox{69.0}{15.0}{10.0}{10.0}{$\scriptstyle{xy=zw}$}
\drawcenteredtext{68.0}{18.0}{$\scriptstyle{x=w}$}
\drawcenteredtext{70.0}{6.0}{$\tau$}
\drawpath{58.0}{18.0}{64.0}{18.0}
\drawpath{58.0}{12.0}{64.0}{12.0}
\drawpath{74.0}{18.0}{80.0}{18.0}
\drawpath{74.0}{12.0}{80.0}{12.0}
\drawcenteredtext{60.0}{20.0}{$x$}
\drawcenteredtext{60.0}{14.0}{$y$}
\drawcenteredtext{78.0}{20.0}{$z$}
\drawcenteredtext{78.0}{14.0}{$w$}
\drawframebox{97.0}{15.0}{10.0}{10.0}{$\scriptstyle{xy=zw}$}
\drawpath{86.0}{18.0}{92.0}{18.0}
\drawpath{86.0}{12.0}{92.0}{12.0}
\drawpath{102.0}{18.0}{108.0}{18.0}
\drawpath{102.0}{12.0}{108.0}{12.0}
\drawcenteredtext{98.0}{18.0}{$\scriptstyle{y=w}$}
\drawcenteredtext{98.0}{6.0}{$\tau^{-1}$}
\drawcenteredtext{44.0}{6.0}{$\eps$}
\drawcenteredtext{86.0}{20.0}{$x$}
\drawcenteredtext{86.0}{10.0}{$y$}
\drawcenteredtext{106.0}{20.0}{$z$}
\drawcenteredtext{106.0}{10.0}{$w$}
\end{picture}
}
\end{center}

Then a colouring of the trefoil, that is, an element of the limit is a tuple of elements of $G$ on the wires satisfying the conditions of the components:

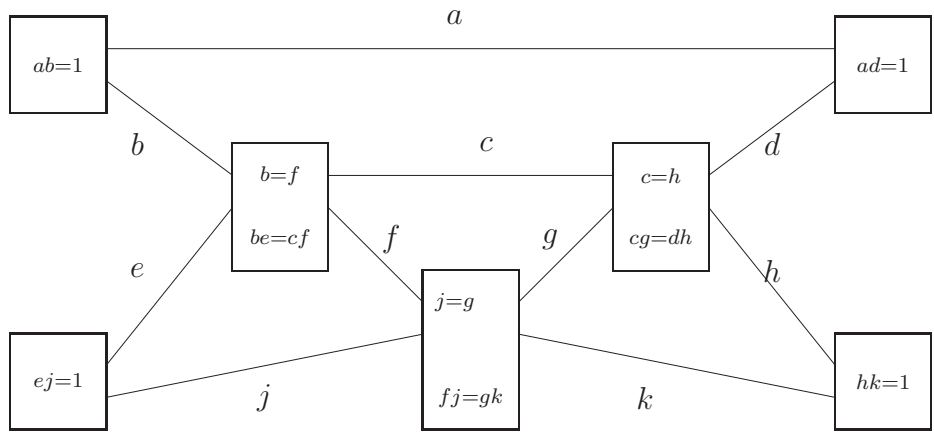
\begin{center}
{\tt\setlength{\unitlength}{6.0pt}
\begin{picture}(66,36)
\thinlines
\drawframebox{7.0}{27.0}{6.0}{6.0}{$\scriptstyle{ab=1}$}
\drawframebox{7.0}{7.0}{6.0}{6.0}{$\scriptstyle{ej=1}$}
\drawframebox{21.0}{18.0}{6.0}{8.0}{}
\drawframebox{33.0}{9.0}{6.0}{10.0}{}
\drawframebox{45.0}{18.0}{6.0}{8.0}{}
\drawframebox{59.0}{27.0}{6.0}{6.0}{$\scriptstyle{ad=1}$}
\drawframebox{59.0}{7.0}{6.0}{6.0}{$\scriptstyle{hk=1}$}
\drawpath{10.0}{28.0}{56.0}{28.0}
\drawpath{10.0}{26.0}{18.0}{20.0}
\drawpath{18.0}{18.0}{10.0}{8.0}
\drawpath{10.0}{6.0}{30.0}{10.0}
\drawpath{24.0}{18.0}{30.0}{12.0}
\drawpath{24.0}{20.0}{42.0}{20.0}
\drawpath{36.0}{12.0}{42.0}{18.0}
\drawpath{56.0}{26.0}{48.0}{20.0}
\drawpath{48.0}{18.0}{56.0}{8.0}
\drawpath{36.0}{10.0}{56.0}{6.0}
\drawcenteredtext{32.0}{30.0}{$a$}
\drawcenteredtext{12.0}{22.0}{$b$}
\drawcenteredtext{34.0}{22.0}{$c$}
\drawcenteredtext{52.0}{22.0}{$d$}
\drawcenteredtext{12.0}{14.0}{$e$}
\drawcenteredtext{28.0}{16.0}{$f$}
\drawcenteredtext{38.0}{16.0}{$g$}
\drawcenteredtext{52.0}{14.0}{$h$}
\drawcenteredtext{20.0}{6.0}{$j$}
\drawcenteredtext{44.0}{6.0}{$k$}
\drawcenteredtext{21.0}{20.0}{$\scriptstyle{b=f}$}
\drawcenteredtext{21.0}{16.0}{$\scriptstyle{be=cf}$}
\drawcenteredtext{45.0}{20.0}{$\scriptstyle{c=h}$}
\drawcenteredtext{45.0}{16.0}{$\scriptstyle{cg=dh}$}
\drawcenteredtext{32.0}{12.0}{$\scriptstyle{j=g}$}
\drawcenteredtext{33.0}{6.0}{$\scriptstyle{fj=gk}$}
\end{picture}
}
\end{center}

When the group is $D_3$ there are 12 colourings, one for each of $(a,c)=(1,1)$, $(123,123)$, $(132,132)$, $(12,12)$,
$(13,13)$, $(23,23)$, $(12,13)$, $(12,23)$, $(13,12)$, $(13,23)$, $(23,12)$, $(23,12)$, whereas the unknot has 6 colourings.

\subsection{Knot groups}\label{subsec-knot-groups}
Consider now the the group object $F$, the free group on one generator, in the category $\gpop$. As we have mentioned the construction $\tsp$ works for any category with finite limits, not just $\set$, and hence there is
a braided monoidal category $\tsp_F(\gpop)$, and a corresponding representation $Gp:\tang\to\tsp_F(\gpop)\to \spa (\gpop)=\csp (\gp)$. 

\begin{thm} (J. Armstrong \cite{aekt})
 If $K$ is a knot then $Gp(K)$ is the knot group of $K$.
\end{thm}
\subsubsection{The knot group of a trefoil}\label{subsubsec-group-trefoil}
\begin{rem}
Limits in $\gpop$ are colimits in $\gp$. We can calculate the knot group from the same picture we used to calculate
the knot colouring. In the diagram

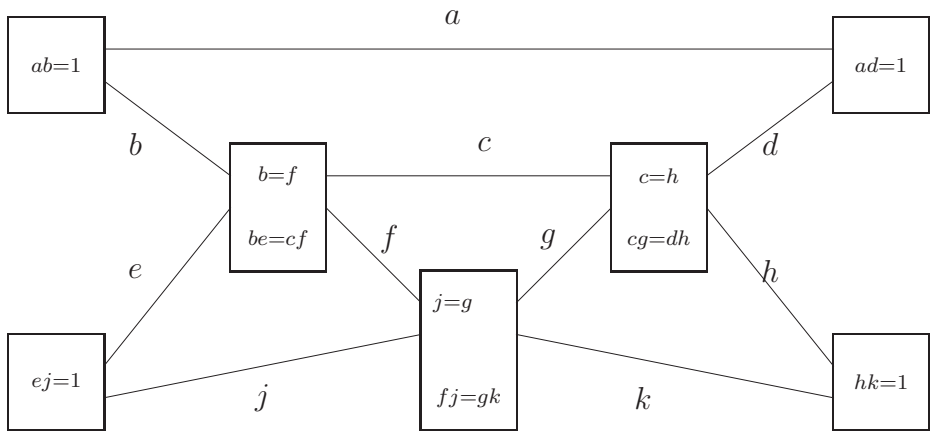
\begin{center}
{\tt\setlength{\unitlength}{6.0pt}
\begin{picture}(66,36)
\thinlines
\drawframebox{7.0}{27.0}{6.0}{6.0}{$\scriptstyle{ab=1}$}
\drawframebox{7.0}{7.0}{6.0}{6.0}{$\scriptstyle{ej=1}$}
\drawframebox{21.0}{18.0}{6.0}{8.0}{}
\drawframebox{33.0}{9.0}{6.0}{10.0}{}
\drawframebox{45.0}{18.0}{6.0}{8.0}{}
\drawframebox{59.0}{27.0}{6.0}{6.0}{$\scriptstyle{ad=1}$}
\drawframebox{59.0}{7.0}{6.0}{6.0}{$\scriptstyle{hk=1}$}
\drawpath{10.0}{28.0}{56.0}{28.0}
\drawpath{10.0}{26.0}{18.0}{20.0}
\drawpath{18.0}{18.0}{10.0}{8.0}
\drawpath{10.0}{6.0}{30.0}{10.0}
\drawpath{24.0}{18.0}{30.0}{12.0}
\drawpath{24.0}{20.0}{42.0}{20.0}
\drawpath{36.0}{12.0}{42.0}{18.0}
\drawpath{56.0}{26.0}{48.0}{20.0}
\drawpath{48.0}{18.0}{56.0}{8.0}
\drawpath{36.0}{10.0}{56.0}{6.0}
\drawcenteredtext{32.0}{30.0}{$a$}
\drawcenteredtext{12.0}{22.0}{$b$}
\drawcenteredtext{34.0}{22.0}{$c$}
\drawcenteredtext{52.0}{22.0}{$d$}
\drawcenteredtext{12.0}{14.0}{$e$}
\drawcenteredtext{28.0}{16.0}{$f$}
\drawcenteredtext{38.0}{16.0}{$g$}
\drawcenteredtext{52.0}{14.0}{$h$}
\drawcenteredtext{20.0}{6.0}{$j$}
\drawcenteredtext{44.0}{6.0}{$k$}
\drawcenteredtext{21.0}{20.0}{$\scriptstyle{b=f}$}
\drawcenteredtext{21.0}{16.0}{$\scriptstyle{be=cf}$}
\drawcenteredtext{45.0}{20.0}{$\scriptstyle{c=h}$}
\drawcenteredtext{45.0}{16.0}{$\scriptstyle{cg=dh}$}
\drawcenteredtext{32.0}{12.0}{$\scriptstyle{j=g}$}
\drawcenteredtext{33.0}{6.0}{$\scriptstyle{fj=gk}$}
\end{picture}
}
\end{center}

\noindent a letter represents the free group $F$ on that generator, letters on a pair of wires represents the free group on two generators $F\times F$ in $\gpop$.
The components are quotients of the free group on the boundary wires by the equations. The evaluation of the circuit
in $\tsp_F(\gpop)$ is a colimit, namely the free group on all the wires quotiented by all the equations.

In the case of the trefoil the knot group is
\begin{align*}
<a,b,c,d,e,f,g,j,k;\;\; & ab=1,b=f,be=cf, c=h,\\
& cg=dh,ad=1,ej=1,j=g,fj=gk,hk=1>.
\end{align*}

\end{rem}

\section{Extending $\tr_G$ and $\tsp_G$}\label{sec-extended}

\subsection{$\tr_{X,G}$}\label{subsec-trel-X,G}
We now describe an extension of $\tr_G$ which depends not only on the group $G$ but also on a set $X$, and we denote it $\tr_{X,G}$, and a similar extension of $\tsp_G$ denoted $\tsp_{X,G}$. These will enable us to model circuits with state.

\begin{definition}
The category $\tr_{X,G}$ has objects $(X\times G)^n$. 
An arrow of $\tr_{X,G}$ is a relation $S$ in $\set$ from 
$ (X\times G)^m$ to $(Y\times G)^n$ 
such that 1) if $(x,h)S(y,k)$ then for any $g\in G$, $(x,h^{g})S(y,k^{g})$, and 2) if $(x,h)S(y,k)$ then $(\ol{h})(\ol{k})^ {-1}\in Z(G)$. Composition and identities are defined as in $\textbf{Rel}$
\end{definition}

In $\tr_{X,G}$ we define a tensor product by $(X\times G)^m \tn (X\times G)^n = (X\times G)^{m+n}$. 

\begin{prop}\label{prop-tcirc-braided}
 $\tr_{X,G}$ is a braided strict monoidal category with 
$$\tw_{(XG)^m \tn (YG)^n}$$ defined to be the relation $$((x,g),(y,h))\sim  
((y,h^{\ol{g}}),(x,g)).$$
\end{prop}

The idea is that in $\tr_{X,G}$ the object $X\times G$ is a single wire carrying data $X$.
As in $\tr_G$ and $\tsp_{G}$, a ``single wire'' $X\times G$ in
$\tr_{X,G}$  admits a commutative Frobenius algebra structure, namely
the comultiplication is the relation $((x,g),(x,h))\sim (x,gh)$; the multiplication is $(x,gh)\sim ((x,g),(x,h))$, the counit is $(x,1)\sim *$ and the unit is $*\sim (x,1)$.

\subsection{Analogue resistive circuits in $\tr_{\bbR,\bbR}$}\label{subsec-resistive}
We begin by describing circuits of resistors which may be described in $\tr_{X,G}$ where $X=\bbR$ is the real numbers, and $G=\bbR$ as a group under addition. 
It is useful to use a graphical notation similar to that of section \ref{subsec-knot-colourings} to do calculations in $\tr_{\bbR,\bbR}$. For example, we draw a relation $S:\bbR\times \bbR\to X\times G$ as:

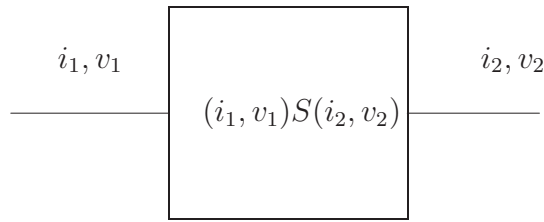
\begin{center}
{\tt\setlength{\unitlength}{5.0pt}
\begin{picture}(52,24)
\thinlines
\drawframebox{25.0}{12.0}{18.0}{16.0}{}
\drawpath{16.0}{12.0}{4.0}{12.0}
\drawpath{34.0}{12.0}{44.0}{12.0}
\drawcenteredtext{10.0}{16.0}{$i_1,v_1$}
\drawcenteredtext{42.0}{16.0}{$i_2,v_2$}
\drawcenteredtext{26.0}{12.0}{$(i_1,v_1) S (i_2,v_2)$}
\end{picture}
}
\end{center}

With this notation, where $i$ denotes current and $v$ denotes voltage, a resistor of resistance $r$ is:

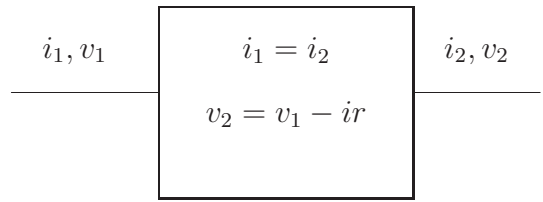
\begin{center}
{\tt\setlength{\unitlength}{4.0pt}
\begin{picture}(58,26)
\thinlines
\drawframebox{30.0}{13.0}{24.0}{18.0}{}
\drawpath{18.0}{14.0}{4.0}{14.0}
\drawpath{42.0}{14.0}{54.0}{14.0}
\drawcenteredtext{10.0}{18.0}{$i_1,v_1$}
\drawcenteredtext{48.0}{18.0}{$i_2,v_2$}
\drawcenteredtext{30.0}{18.0}{$i_1=i_2$}
\drawcenteredtext{30.0}{12.0}{$v_2=v_1-ir$}
\end{picture}
}
\end{center}

The unit and counit, which sometimes we draw as forks, and which embody Kirchhoff's law of currents:

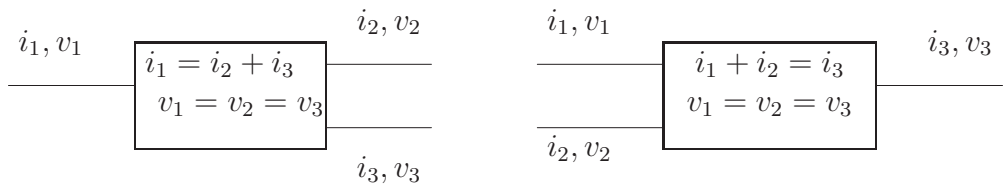
\begin{center}
{\tt\setlength{\unitlength}{4pt} 
\begin{picture}(104,26)
\thinlines
\drawframebox{25.0}{13.0}{18.0}{10.0}{}
\drawpath{16.0}{14.0}{10.0}{14.0}
\drawpath{10.0}{14.0}{4.0}{14.0}
\drawpath{34.0}{16.0}{44.0}{16.0}
\drawpath{34.0}{10.0}{44.0}{10.0}
\drawframebox{76.0}{13.0}{20.0}{10.0}{}
\drawpath{66.0}{16.0}{54.0}{16.0}
\drawpath{66.0}{10.0}{54.0}{10.0}
\drawpath{86.0}{14.0}{98.0}{14.0}
\drawcenteredtext{8.0}{18.0}{$i_1,v_1$}
\drawcenteredtext{40.0}{20.0}{$i_2,v_2$}
\drawcenteredtext{40.0}{6.0}{$i_3,v_3$}
\drawcenteredtext{58.0}{20.0}{$i_1,v_1$}
\drawcenteredtext{58.0}{8.0}{$i_2,v_2$}
\drawcenteredtext{94.0}{18.0}{$i_3,v_3$}
\drawcenteredtext{24.0}{16.0}{$i_1=i_2+i_3$}
\drawcenteredtext{26.0}{12.0}{$v_1=v_2=v_3$}
\drawcenteredtext{76.0}{16.0}{$i_1+i_2=i_3$}
\drawcenteredtext{76.0}{12.0}{$v_1=v_2=v_3$}
\end{picture}
}
\end{center}

Using the operations of $\tr_{\bbR,\bbR}$ one can now evaluate a network of resistors. For example the circuit with two parallel resistors with resistances $r_1$, $r_2$ respectively

\begin{center}
{\tt\setlength{\unitlength}{3.0pt}
\begin{picture}(56,30)
\thinlines
\drawpath{4.0}{12.0}{14.0}{12.0}
\drawpath{14.0}{12.0}{20.0}{18.0}
\drawpath{20.0}{18.0}{14.0}{12.0}
\drawpath{14.0}{12.0}{20.0}{6.0}
\drawpath{20.0}{6.0}{22.0}{6.0}
\drawpath{22.0}{6.0}{24.0}{8.0}
\drawpath{24.0}{8.0}{24.0}{4.0}
\drawpath{24.0}{4.0}{26.0}{8.0}
\drawpath{26.0}{8.0}{26.0}{4.0}
\drawpath{26.0}{4.0}{28.0}{8.0}
\drawpath{28.0}{8.0}{28.0}{4.0}
\drawpath{28.0}{4.0}{30.0}{8.0}
\drawpath{30.0}{8.0}{30.0}{4.0}
\drawpath{30.0}{4.0}{32.0}{8.0}
\drawpath{32.0}{8.0}{32.0}{4.0}
\drawpath{32.0}{4.0}{34.0}{8.0}
\drawpath{34.0}{8.0}{34.0}{4.0}
\drawpath{34.0}{4.0}{36.0}{6.0}
\drawpath{36.0}{6.0}{38.0}{6.0}
\drawpath{38.0}{6.0}{44.0}{12.0}
\drawpath{44.0}{12.0}{38.0}{18.0}
\drawpath{38.0}{18.0}{38.0}{18.0}
\drawpath{44.0}{12.0}{52.0}{12.0}
\drawpath{20.0}{18.0}{24.0}{18.0}
\drawpath{24.0}{18.0}{26.0}{20.0}
\drawpath{26.0}{20.0}{26.0}{16.0}
\drawpath{26.0}{16.0}{28.0}{20.0}
\drawpath{28.0}{20.0}{28.0}{16.0}
\drawpath{28.0}{16.0}{30.0}{20.0}
\drawpath{30.0}{20.0}{30.0}{16.0}
\drawpath{30.0}{16.0}{32.0}{20.0}
\drawpath{32.0}{20.0}{32.0}{16.0}
\drawpath{32.0}{16.0}{34.0}{20.0}
\drawpath{34.0}{20.0}{34.0}{16.0}
\drawpath{34.0}{16.0}{36.0}{18.0}
\drawpath{36.0}{18.0}{38.0}{18.0}
\drawpath{36.0}{16.0}{36.0}{16.0}
\drawcenteredtext{30.0}{24.0}{$r_1$}
\drawcenteredtext{30.0}{10.0}{$r_2$}
\end{picture}
}
\end{center}

 evaluates as:

\begin{center}
 {\tt\setlength{\unitlength}{5.0pt}
\begin{picture}(50,24)
\thinlines
\drawframebox{24.0}{12.0}{20.0}{16.0}{}
\drawpath{14.0}{12.0}{4.0}{12.0}
\drawpath{34.0}{12.0}{44.0}{12.0}
\drawcenteredtext{8.0}{14.0}{$i_1,v_1$}
\drawcenteredtext{40.0}{14.0}{$i_2,v_2$}
\drawcenteredtext{24.0}{18.0}{$i_1=i_2$}
\drawcenteredtext{25.0}{10.0}{$v_2-v_1=i_1({\frac {r_1r_2} {r_1+r_2}})$}
\end{picture}
}
\end{center}

\subsection{$\tsp_{X,G}$}\label{subsec-tspan-XG}
\begin{definition}
The category $\tsp_{X,G}$ has objects $(X\times G)^n$. 
An arrow of $\tsp_{X,G}$ is an isomorphism of spans $S$ in $\set$:
$ (X\times G)^m \toleft^{\delta_0} S \to^{\delta_1} (Y\times G)^n$ 
such that
such that there exist a function $G\times S \to S$ of $G$ on $S$ written $(g,s)\mapsto gs$ yielding a bijection for each $g\in G$, and 
satisfying:
\begin{itemize}
\item[1)] if $\delta_0(s) = (x_1,h_1,...,x_m,h_m)$ and $\delta_1(s) = (y_1,k_1,...,y_n,k_n)$ then %
$\delta_0(gs) = (x_1,h_1^g,...,x_m,h_m^g)$  and $\delta_1(gs) = (y_1,k_1^g, ... ,y_n,k_n^g)$ for all $g$ in $G$,
\item[2)] if $\delta_0(s) = (x_1,h_1,...,x_m,h_m)$ and $\delta_1(s) = (y_1,k_1,...y_n,k_n)$  then $h_1...h_m(k_1...,k_n)^{-1}\in Z(G)$.
\end{itemize}
 Composition and identities and tensor are defined as in $\textbf{Span}$.
The braiding and Frobenius structure are as in $\tr_{X,G}$.  
\end{definition}

It is clear that this definition may be made in any category $\bC$ with finite limits
to give a category $\tsp_{X,G}(\bC)$.

\subsection{RLC circuits in $\tsp_{X,G}(\gra)$}\label{subsec-RLC}
The algebra of RLC circuits we will describe was introduced in \cite{KSW00a} but without the consciousness of Frobenius algebras. We will give a brief recapitulation without full details.

We need to say something first about the somewhat unusual interpretation of a graph in this setting. If the graph consists of the two (domain and codomain) functions $\phi: X\to Y$ and $\psi: X\to Y$ we will interpret this as the formal differential equation $\phi^{'} =\psi$. For further explanation of this interpretation see \cite{KSW00a}. In the examples we describe the interpretation will have a clear meaning. There is a notion of behaviour for such a system, namely a function $x:\bbR \to X$ such that $\phi^{'}(x(t))=\psi(x(t))$ (only meaningful with smoothness assumptions).

We will now consider $\tsp_{X,G}(\gra)$ where both $X$ and $G$ are the graph with one vertex, and set of arrows $\bbR$; we will identify both $X$ and $G$ with the set $\bbR$, the group structure being addition. 

Again it is useful to use a graphical notation similar to that of Section \ref{subsec-knot-colourings} to do calculations in $\tsp_{\bbR,\bbR}$. For example, we draw the spans corresponding to the constants $\tw$, $\Delta$, $\nabla$, $\eta$, $\eps$, and the resistors of the algebra (in which all of the graphs have one vertex) exactly as in section \ref{subsec-resistive}. 

Instead the graph of a capacitor with capacitance $c$ is the pair functions $\phi,\psi:\bbR^3\to \bbR$ defined by
$\phi(i,v,q)=q$ and $\psi (i,v,q)=i$; the interpretation of this is that a capacitor has state $i$, $v$, and also state $q$, the charge of the capacitor, and that $q^{'}=i$. The boundary conditions (the morphism of the span)
are on the left $v_1=v$, and  $i_1=i$ and on the right $v_2=v-{\frac q c}$ and $i_2=i$. Hence we draw the capacitor as follows:

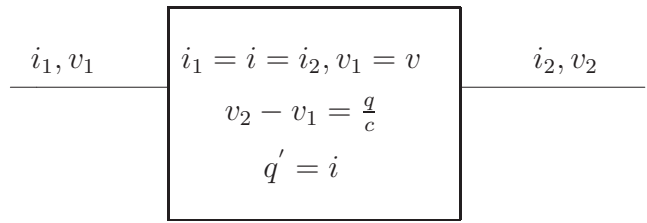
\begin{center}
{\tt\setlength{\unitlength}{5.0pt}
\begin{picture}(56,24)
\thinlines
\drawcenteredtext{26.0}{8.0}{$q^{'}=i$}
\drawpath{16.0}{14.0}{6.0}{14.0}
\drawpath{38.0}{14.0}{52.0}{14.0}
\drawpath{6.0}{14.0}{4.0}{14.0}
\drawcenteredtext{8.0}{16.0}{$i_1,v_1$}
\drawcenteredtext{46.0}{16.0}{$i_2,v_2$}
\drawcenteredtext{26.0}{16.0}{$i_1=i=i_2, v_1=v$}
\drawcenteredtext{26.0}{12.0}{$v_2-v_1={\frac {q}{c}}$}
\drawframebox{27.0}{12.0}{22.0}{16.0}{}
\end{picture}
}
\end{center}

 Similarly an inductor with inductance $l$ has an extra variable of state $p$ with graph $\bbR^3\to \bbR$, and pictures;

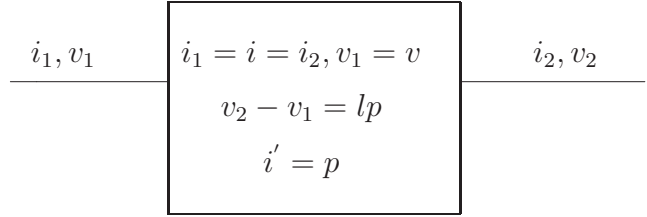
\begin{center}
{\tt\setlength{\unitlength}{5.0pt}
\begin{picture}(56,24)
\thinlines
\drawcenteredtext{26.0}{8.0}{$i^{'}=p$}
\drawpath{16.0}{14.0}{6.0}{14.0}
\drawpath{38.0}{14.0}{52.0}{14.0}
\drawpath{6.0}{14.0}{4.0}{14.0}
\drawcenteredtext{8.0}{16.0}{$i_1,v_1$}
\drawcenteredtext{46.0}{16.0}{$i_2,v_2$}
\drawcenteredtext{26.0}{16.0}{$i_1=i=i_2, v_1=v$}
\drawcenteredtext{26.0}{12.0}{$v_2-v_1=lp$}
\drawframebox{27.0}{12.0}{22.0}{16.0}{}
\end{picture}
}
\end{center}
 
 Using the operations of $\tsp_{\bbR,\bbR}(\gra)$ one can now evaluate a network of resistors, capacitors and inductors. For example the circuit of an inductance and a capacitance
 
 \begin{center}
 {\tt\setlength{\unitlength}{4.0pt}
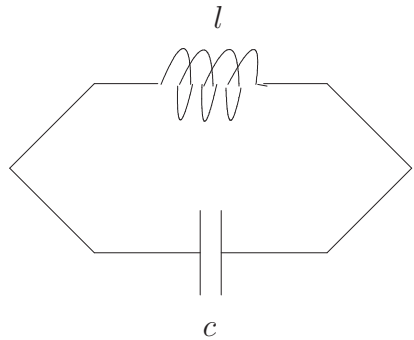
\begin{picture}(46,42)
\thinlines
\drawpath{12.0}{30.0}{4.0}{22.0}
\drawpath{4.0}{22.0}{12.0}{14.0}
\drawpath{12.0}{14.0}{22.0}{14.0}
\drawpath{12.0}{30.0}{18.0}{30.0}
\drawpath{22.0}{18.0}{22.0}{10.0}
\drawpath{24.0}{18.0}{24.0}{10.0}
\drawpath{24.0}{14.0}{34.0}{14.0}
\drawpath{34.0}{14.0}{42.0}{22.0}
\drawpath{42.0}{22.0}{34.0}{30.0}
\drawpath{34.0}{30.0}{28.0}{30.0}
\path(18.31,29.93)(18.31,29.93)(18.36,30.06)(18.41,30.19)(18.47,30.33)(18.52,30.45)(18.58,30.58)(18.63,30.69)(18.68,30.81)(18.74,30.93)
\path(18.74,30.93)(18.79,31.05)(18.84,31.16)(18.88,31.26)(18.93,31.37)(18.99,31.47)(19.02,31.56)(19.08,31.66)(19.13,31.76)(19.18,31.84)
\path(19.18,31.84)(19.22,31.94)(19.27,32.01)(19.31,32.11)(19.36,32.19)(19.4,32.26)(19.45,32.33)(19.49,32.41)(19.52,32.48)(19.56,32.55)
\path(19.56,32.55)(19.61,32.61)(19.65,32.66)(19.68,32.73)(19.72,32.79)(19.77,32.83)(19.81,32.88)(19.84,32.94)(19.88,32.98)(19.91,33.01)
\path(19.91,33.01)(19.95,33.05)(20.0,33.09)(20.02,33.12)(20.06,33.16)(20.09,33.19)(20.13,33.22)(20.15,33.23)(20.2,33.26)(20.22,33.27)
\path(20.22,33.27)(20.25,33.3)(20.29,33.3)(20.31,33.31)(20.34,33.31)(20.38,33.33)(20.4,33.33)(20.43,33.33)(20.45,33.31)(20.49,33.31)
\path(20.49,33.31)(20.52,33.3)(20.54,33.3)(20.56,33.27)(20.59,33.26)(20.61,33.23)(20.63,33.22)(20.65,33.19)(20.68,33.16)(20.7,33.12)
\path(20.7,33.12)(20.72,33.09)(20.75,33.05)(20.77,33.01)(20.79,32.98)(20.79,32.94)(20.81,32.88)(20.84,32.83)(20.86,32.79)(20.86,32.73)
\path(20.86,32.73)(20.88,32.66)(20.9,32.61)(20.91,32.55)(20.93,32.48)(20.95,32.41)(20.95,32.33)(20.97,32.26)(20.99,32.19)(21.0,32.11)
\path(21.0,32.11)(21.0,32.01)(21.02,31.94)(21.02,31.84)(21.04,31.76)(21.04,31.66)(21.06,31.56)(21.06,31.47)(21.06,31.37)(21.08,31.26)
\path(21.08,31.26)(21.09,31.16)(21.09,31.05)(21.09,30.93)(21.09,30.81)(21.09,30.69)(21.11,30.58)(21.11,30.45)(21.11,30.33)(21.11,30.19)
\path(21.11,30.19)(21.11,30.06)(21.11,29.93)(21.11,29.93)
\path(21.29,29.93)(21.29,29.93)(21.25,29.79)(21.22,29.68)(21.2,29.54)(21.16,29.43)(21.13,29.29)(21.11,29.18)(21.08,29.06)(21.04,28.95)
\path(21.04,28.95)(21.02,28.84)(21.0,28.74)(20.97,28.63)(20.93,28.52)(20.9,28.43)(20.88,28.34)(20.86,28.24)(20.83,28.15)(20.81,28.06)
\path(20.81,28.06)(20.77,27.97)(20.75,27.88)(20.72,27.81)(20.7,27.72)(20.68,27.65)(20.65,27.58)(20.63,27.5)(20.61,27.43)(20.59,27.36)
\path(20.59,27.36)(20.56,27.31)(20.54,27.25)(20.52,27.18)(20.5,27.13)(20.47,27.08)(20.45,27.02)(20.43,26.97)(20.41,26.93)(20.4,26.9)
\path(20.4,26.9)(20.38,26.86)(20.36,26.81)(20.34,26.79)(20.33,26.75)(20.31,26.72)(20.29,26.7)(20.27,26.66)(20.25,26.65)(20.24,26.63)
\path(20.24,26.63)(20.22,26.61)(20.2,26.61)(20.2,26.59)(20.18,26.59)(20.16,26.58)(20.15,26.58)(20.13,26.58)(20.13,26.58)(20.11,26.58)
\path(20.11,26.58)(20.09,26.59)(20.09,26.59)(20.06,26.61)(20.06,26.63)(20.04,26.65)(20.04,26.66)(20.02,26.68)(20.02,26.72)(20.0,26.75)
\path(20.0,26.75)(20.0,26.77)(19.99,26.81)(19.97,26.84)(19.97,26.88)(19.95,26.93)(19.95,26.97)(19.95,27.02)(19.93,27.06)(19.93,27.11)
\path(19.93,27.11)(19.93,27.16)(19.91,27.22)(19.9,27.29)(19.9,27.34)(19.9,27.41)(19.9,27.49)(19.88,27.56)(19.88,27.63)(19.88,27.7)
\path(19.88,27.7)(19.88,27.77)(19.88,27.86)(19.88,27.95)(19.88,28.02)(19.86,28.11)(19.86,28.2)(19.86,28.31)(19.86,28.4)(19.86,28.5)
\path(19.86,28.5)(19.86,28.59)(19.86,28.7)(19.86,28.81)(19.86,28.91)(19.86,29.04)(19.86,29.15)(19.86,29.27)(19.88,29.38)(19.88,29.5)
\path(19.88,29.5)(19.88,29.63)(19.88,29.75)(19.88,29.77)
\path(20.06,29.93)(20.06,29.93)(20.11,30.06)(20.18,30.19)(20.24,30.31)(20.29,30.44)(20.36,30.55)(20.4,30.68)(20.47,30.79)(20.52,30.91)
\path(20.52,30.91)(20.58,31.01)(20.63,31.12)(20.68,31.23)(20.74,31.33)(20.79,31.43)(20.84,31.51)(20.9,31.62)(20.95,31.7)(21.0,31.8)
\path(21.0,31.8)(21.04,31.87)(21.09,31.95)(21.15,32.04)(21.2,32.12)(21.24,32.19)(21.29,32.26)(21.34,32.33)(21.38,32.4)(21.43,32.47)
\path(21.43,32.47)(21.47,32.52)(21.52,32.58)(21.56,32.63)(21.61,32.69)(21.65,32.75)(21.68,32.79)(21.72,32.83)(21.77,32.87)(21.81,32.91)
\path(21.81,32.91)(21.84,32.95)(21.88,32.98)(21.93,33.02)(21.97,33.05)(22.0,33.08)(22.04,33.09)(22.08,33.12)(22.11,33.13)(22.15,33.16)
\path(22.15,33.16)(22.18,33.16)(22.22,33.18)(22.25,33.19)(22.29,33.19)(22.31,33.19)(22.34,33.19)(22.38,33.19)(22.4,33.19)(22.43,33.18)
\path(22.43,33.18)(22.47,33.16)(22.5,33.16)(22.52,33.13)(22.56,33.12)(22.58,33.09)(22.61,33.06)(22.63,33.05)(22.65,33.01)(22.68,32.98)
\path(22.68,32.98)(22.7,32.94)(22.74,32.91)(22.75,32.87)(22.77,32.83)(22.81,32.77)(22.83,32.73)(22.84,32.68)(22.86,32.62)(22.88,32.56)
\path(22.88,32.56)(22.9,32.51)(22.93,32.44)(22.95,32.38)(22.95,32.31)(22.97,32.25)(23.0,32.18)(23.0,32.09)(23.02,32.01)(23.04,31.94)
\path(23.04,31.94)(23.06,31.86)(23.06,31.76)(23.09,31.68)(23.09,31.58)(23.11,31.5)(23.11,31.4)(23.13,31.3)(23.13,31.19)(23.15,31.08)
\path(23.15,31.08)(23.15,30.98)(23.16,30.87)(23.18,30.76)(23.18,30.65)(23.18,30.52)(23.2,30.41)(23.2,30.29)(23.2,30.16)(23.22,30.02)
\path(23.22,30.02)(23.22,29.9)(23.22,29.77)(23.22,29.77)
\path(23.58,29.93)(23.58,29.93)(23.54,29.79)(23.5,29.65)(23.47,29.52)(23.45,29.4)(23.41,29.27)(23.38,29.15)(23.36,29.02)(23.33,28.9)
\path(23.33,28.9)(23.29,28.79)(23.27,28.68)(23.25,28.56)(23.22,28.45)(23.18,28.36)(23.16,28.25)(23.13,28.15)(23.11,28.06)(23.09,27.97)
\path(23.09,27.97)(23.06,27.88)(23.04,27.79)(23.0,27.7)(22.99,27.61)(22.95,27.54)(22.93,27.47)(22.91,27.38)(22.9,27.31)(22.86,27.25)
\path(22.86,27.25)(22.84,27.18)(22.83,27.11)(22.81,27.06)(22.79,27.0)(22.75,26.95)(22.74,26.9)(22.72,26.84)(22.7,26.81)(22.68,26.75)
\path(22.68,26.75)(22.66,26.72)(22.65,26.68)(22.63,26.65)(22.61,26.61)(22.59,26.59)(22.58,26.56)(22.56,26.54)(22.54,26.52)(22.52,26.5)
\path(22.52,26.5)(22.5,26.47)(22.49,26.47)(22.47,26.45)(22.45,26.45)(22.45,26.45)(22.43,26.45)(22.41,26.45)(22.4,26.45)(22.38,26.45)
\path(22.38,26.45)(22.38,26.47)(22.36,26.47)(22.36,26.5)(22.34,26.52)(22.33,26.54)(22.31,26.56)(22.31,26.59)(22.29,26.61)(22.29,26.65)
\path(22.29,26.65)(22.27,26.68)(22.27,26.72)(22.25,26.75)(22.25,26.81)(22.25,26.84)(22.24,26.9)(22.22,26.95)(22.22,27.0)(22.22,27.06)
\path(22.22,27.06)(22.2,27.11)(22.2,27.18)(22.2,27.25)(22.2,27.31)(22.18,27.38)(22.18,27.47)(22.18,27.54)(22.18,27.61)(22.16,27.7)
\path(22.16,27.7)(22.16,27.79)(22.16,27.88)(22.15,27.97)(22.15,28.06)(22.15,28.15)(22.15,28.25)(22.15,28.36)(22.15,28.45)(22.15,28.56)
\path(22.15,28.56)(22.15,28.68)(22.15,28.79)(22.15,28.9)(22.15,29.02)(22.15,29.15)(22.15,29.27)(22.15,29.4)(22.16,29.52)(22.16,29.65)
\path(22.16,29.65)(22.16,29.79)(22.18,29.93)(22.18,29.93)
\path(22.34,29.77)(22.34,29.77)(22.4,29.9)(22.47,30.04)(22.52,30.16)(22.59,30.3)(22.65,30.43)(22.7,30.55)(22.77,30.66)(22.83,30.79)
\path(22.83,30.79)(22.88,30.91)(22.95,31.01)(23.0,31.12)(23.06,31.23)(23.11,31.33)(23.16,31.44)(23.22,31.54)(23.27,31.63)(23.34,31.73)
\path(23.34,31.73)(23.38,31.81)(23.43,31.91)(23.49,31.98)(23.54,32.06)(23.59,32.15)(23.63,32.23)(23.68,32.3)(23.74,32.37)(23.79,32.44)
\path(23.79,32.44)(23.84,32.51)(23.88,32.56)(23.93,32.62)(23.97,32.69)(24.02,32.73)(24.06,32.79)(24.11,32.83)(24.15,32.88)(24.2,32.93)
\path(24.2,32.93)(24.24,32.97)(24.27,33.01)(24.31,33.04)(24.36,33.08)(24.4,33.11)(24.43,33.12)(24.47,33.16)(24.52,33.18)(24.54,33.19)
\path(24.54,33.19)(24.59,33.2)(24.63,33.22)(24.65,33.23)(24.7,33.23)(24.72,33.23)(24.75,33.25)(24.79,33.23)(24.83,33.23)(24.86,33.23)
\path(24.86,33.23)(24.88,33.22)(24.91,33.2)(24.95,33.19)(24.97,33.18)(25.0,33.16)(25.04,33.12)(25.06,33.11)(25.09,33.08)(25.11,33.04)
\path(25.11,33.04)(25.15,33.01)(25.16,32.97)(25.2,32.93)(25.22,32.88)(25.25,32.83)(25.27,32.79)(25.29,32.73)(25.31,32.69)(25.33,32.62)
\path(25.33,32.62)(25.34,32.56)(25.36,32.51)(25.38,32.44)(25.4,32.37)(25.43,32.3)(25.45,32.23)(25.45,32.15)(25.47,32.06)(25.5,31.98)
\path(25.5,31.98)(25.5,31.91)(25.52,31.81)(25.54,31.73)(25.54,31.63)(25.56,31.54)(25.58,31.44)(25.59,31.33)(25.59,31.23)(25.61,31.12)
\path(25.61,31.12)(25.61,31.01)(25.63,30.91)(25.63,30.79)(25.65,30.66)(25.65,30.55)(25.65,30.43)(25.66,30.3)(25.66,30.16)(25.68,30.04)
\path(25.68,30.04)(25.68,29.9)(25.68,29.77)(25.68,29.77)
\path(25.86,29.59)(25.86,29.59)(25.84,29.45)(25.81,29.34)(25.79,29.22)(25.77,29.11)(25.74,29.0)(25.72,28.9)(25.7,28.79)(25.66,28.68)
\path(25.66,28.68)(25.65,28.59)(25.63,28.49)(25.61,28.38)(25.58,28.29)(25.56,28.2)(25.54,28.11)(25.52,28.02)(25.5,27.95)(25.47,27.86)
\path(25.47,27.86)(25.45,27.79)(25.43,27.7)(25.4,27.63)(25.38,27.56)(25.36,27.49)(25.34,27.43)(25.33,27.36)(25.31,27.29)(25.29,27.24)
\path(25.29,27.24)(25.27,27.18)(25.25,27.13)(25.24,27.08)(25.22,27.02)(25.2,26.97)(25.18,26.93)(25.15,26.88)(25.15,26.84)(25.13,26.81)
\path(25.13,26.81)(25.11,26.77)(25.09,26.75)(25.08,26.72)(25.06,26.68)(25.04,26.65)(25.02,26.63)(25.0,26.61)(25.0,26.59)(24.97,26.58)
\path(24.97,26.58)(24.97,26.56)(24.95,26.56)(24.93,26.54)(24.91,26.54)(24.9,26.54)(24.88,26.54)(24.88,26.54)(24.86,26.54)(24.84,26.54)
\path(24.84,26.54)(24.84,26.56)(24.81,26.56)(24.81,26.58)(24.79,26.59)(24.79,26.61)(24.77,26.63)(24.75,26.65)(24.75,26.68)(24.74,26.72)
\path(24.74,26.72)(24.72,26.75)(24.72,26.77)(24.7,26.81)(24.68,26.84)(24.68,26.88)(24.66,26.93)(24.65,26.97)(24.65,27.02)(24.63,27.08)
\path(24.63,27.08)(24.63,27.13)(24.61,27.18)(24.61,27.24)(24.61,27.29)(24.59,27.36)(24.59,27.43)(24.58,27.49)(24.56,27.56)(24.56,27.63)
\path(24.56,27.63)(24.56,27.7)(24.54,27.79)(24.54,27.86)(24.54,27.95)(24.52,28.02)(24.52,28.11)(24.52,28.2)(24.5,28.29)(24.5,28.38)
\path(24.5,28.38)(24.5,28.49)(24.49,28.59)(24.49,28.68)(24.47,28.79)(24.47,28.9)(24.47,29.0)(24.47,29.11)(24.47,29.22)(24.45,29.34)
\path(24.45,29.34)(24.45,29.45)(24.45,29.58)(24.45,29.59)
\path(24.63,29.93)(24.63,29.93)(24.7,30.06)(24.77,30.19)(24.84,30.31)(24.9,30.44)(24.97,30.55)(25.04,30.68)(25.11,30.8)(25.16,30.91)
\path(25.16,30.91)(25.22,31.01)(25.29,31.12)(25.36,31.23)(25.41,31.33)(25.47,31.43)(25.52,31.52)(25.59,31.62)(25.65,31.7)(25.7,31.8)
\path(25.7,31.8)(25.75,31.88)(25.81,31.97)(25.86,32.05)(25.91,32.12)(25.97,32.19)(26.02,32.27)(26.06,32.34)(26.11,32.41)(26.16,32.48)
\path(26.16,32.48)(26.22,32.54)(26.25,32.59)(26.31,32.66)(26.34,32.7)(26.38,32.76)(26.43,32.8)(26.47,32.86)(26.52,32.9)(26.56,32.94)
\path(26.56,32.94)(26.59,32.98)(26.63,33.01)(26.66,33.05)(26.7,33.08)(26.75,33.11)(26.77,33.12)(26.81,33.16)(26.84,33.18)(26.88,33.19)
\path(26.88,33.19)(26.9,33.2)(26.93,33.22)(26.97,33.23)(27.0,33.23)(27.02,33.23)(27.06,33.23)(27.08,33.23)(27.11,33.23)(27.13,33.23)
\path(27.13,33.23)(27.15,33.22)(27.18,33.2)(27.2,33.19)(27.22,33.18)(27.24,33.16)(27.25,33.12)(27.27,33.11)(27.29,33.08)(27.31,33.05)
\path(27.31,33.05)(27.31,33.01)(27.34,32.98)(27.34,32.94)(27.36,32.9)(27.38,32.86)(27.38,32.8)(27.4,32.76)(27.4,32.7)(27.41,32.66)
\path(27.41,32.66)(27.41,32.59)(27.43,32.54)(27.43,32.48)(27.43,32.41)(27.43,32.34)(27.45,32.27)(27.45,32.19)(27.45,32.12)(27.45,32.05)
\path(27.45,32.05)(27.45,31.97)(27.45,31.88)(27.45,31.8)(27.43,31.7)(27.43,31.62)(27.43,31.52)(27.43,31.43)(27.41,31.33)(27.41,31.23)
\path(27.41,31.23)(27.4,31.12)(27.4,31.01)(27.38,30.91)(27.38,30.8)(27.36,30.68)(27.34,30.55)(27.34,30.44)(27.31,30.31)(27.31,30.19)
\path(27.31,30.19)(27.29,30.06)(27.27,29.93)(27.27,29.93)
\drawpath{27.45}{29.93}{28.33}{29.77}
\drawcenteredtext{23.75}{36.38}{$l$}
\drawcenteredtext{22.88}{6.77}{$c$}
\end{picture}
}
 \end{center}

evaluates as 

\begin{center}
{\tt\setlength{\unitlength}{4.0pt}
\begin{picture}(36,31)
\thinlines
\drawframebox{18.89}{15.96}{28.12}{22.12}{}
\drawcenteredtext{17.65}{23.72}{$i,v_1,v_2,p,q$}
\drawcenteredtext{17.13}{18.15}{$(-i)'=p$}
\drawcenteredtext{17.13}{13.09}{$q'=i$}
\drawcenteredtext{20}{7.86}{${\frac q c}=v_1-v_2=lp$}
\end{picture}
}
\end{center}

A behaviour consists of five functions from $\bbR$ to $\bbR$, namely $i(t)$, $v_1(t)$, $v_2(t)$,  $q(t)$, $p(t)$ such that $i^{'}=-p$, $q^{'}=i$ and ${\frac q c}=v_2-v_1=lp$. 

\section{Dirac's belt trick}
The claim is that the following two circuits are equal in $\tcd$, that is that a rotation through $2\pi$ of a component $I\to X^3$ is equal to the identity. We suspect but are unable to prove that a rotation through $\pi$ is not the identity - however in $\tr_G$ it is.
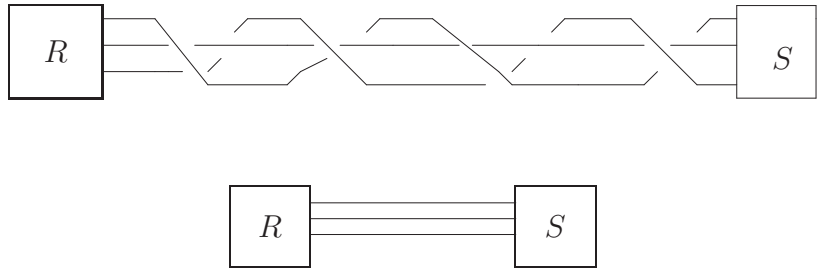
\begin{center}
{\tt\setlength{\unitlength}{2.5pt}
\begin{picture}(130,22)
\thinlines
\drawframebox{11.0}{11.0}{14.0}{14.0}{$R$}
\drawpath{18.0}{16.0}{26.0}{16.0}
\drawpath{26.0}{16.0}{34.0}{6.0}
\drawpath{34.0}{6.0}{46.0}{6.0}
\drawpath{46.0}{6.0}{48.0}{8.0}
\drawpath{18.0}{12.0}{28.0}{12.0}
\drawpath{32.0}{12.0}{46.0}{12.0}
\drawpath{18.0}{8.0}{26.0}{8.0}
\drawpath{26.0}{8.0}{28.0}{8.0}
\drawpath{28.0}{8.0}{30.0}{8.0}
\drawpath{30.0}{8.0}{30.0}{8.0}
\drawpath{34.0}{8.0}{36.0}{10.0}
\drawpath{38.0}{14.0}{40.0}{16.0}
\drawpath{40.0}{16.0}{48.0}{16.0}
\drawpath{48.0}{16.0}{58.0}{6.0}
\drawpath{46.0}{12.0}{50.0}{12.0}
\drawpath{54.0}{12.0}{62.0}{12.0}
\drawpath{48.0}{8.0}{52.0}{10.0}
\drawpath{58.0}{14.0}{60.0}{16.0}
\drawpath{60.0}{16.0}{68.0}{16.0}
\drawpath{68.0}{16.0}{78.0}{8.0}
\drawpath{62.0}{12.0}{72.0}{12.0}
\drawpath{74.0}{12.0}{84.0}{12.0}
\drawpath{58.0}{6.0}{76.0}{6.0}
\drawpath{78.0}{8.0}{80.0}{6.0}
\drawpath{80.0}{6.0}{90.0}{6.0}
\drawpath{80.0}{8.0}{82.0}{10.0}
\drawpath{86.0}{14.0}{88.0}{16.0}
\drawpath{88.0}{16.0}{98.0}{16.0}
\drawpath{98.0}{16.0}{108.0}{6.0}
\drawpath{84.0}{12.0}{100.0}{12.0}
\drawpath{104.0}{12.0}{114.0}{12.0}
\drawpath{90.0}{6.0}{100.0}{6.0}
\drawpath{100.0}{6.0}{102.0}{8.0}
\drawpath{108.0}{14.0}{110.0}{16.0}
\drawpath{110.0}{16.0}{114.0}{16.0}
\drawpath{108.0}{6.0}{114.0}{6.0}
\drawpath{114.0}{18.0}{114.0}{4.0}
\drawpath{114.0}{4.0}{126.0}{4.0}
\drawpath{126.0}{4.0}{126.0}{16.0}
\drawpath{126.0}{16.0}{126.0}{18.0}
\drawpath{126.0}{18.0}{114.0}{18.0}
\drawcenteredtext{121.0}{10.0}{$S$}
\end{picture}
}
\end{center}

\begin{center}
{\tt\setlength{\unitlength}{3.0pt}
\begin{picture}(54,18)
\thinlines
\drawframebox{9.0}{9.0}{10.0}{10.0}{$R$}
\drawframebox{45.0}{9.0}{10.0}{10.0}{$S$}
\drawpath{14.0}{12.0}{40.0}{12.0}
\drawpath{14.0}{10.0}{40.0}{10.0}
\drawpath{14.0}{8.0}{40.0}{8.0}
\drawcenteredtext{10.0}{8.0}{}
\end{picture}
}
\end{center}

We  give a sketch of a proof only.
Using arguments similar to that of example \ref{subsubsec-flash} we may prove that the first (twisted) circuit is equal to

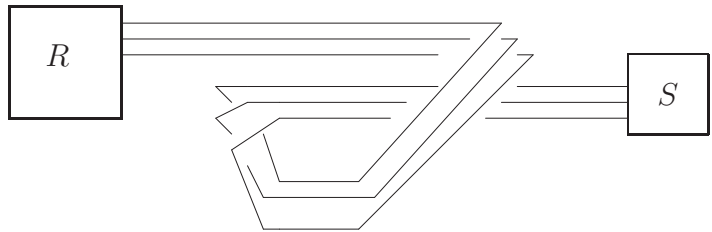
\begin{center}
{\tt\setlength{\unitlength}{3.0pt}
\begin{picture}(96,36)
\thinlines
\drawframebox{11.0}{25.0}{14.0}{14.0}{}
\drawpath{18.0}{30.0}{66.0}{30.0}
\drawpath{66.0}{30.0}{48.0}{10.0}
\drawpath{48.0}{10.0}{38.0}{10.0}
\drawpath{38.0}{10.0}{36.0}{16.0}
\drawpath{18.0}{28.0}{62.0}{28.0}
\drawpath{66.0}{28.0}{68.0}{28.0}
\drawpath{68.0}{28.0}{50.0}{8.0}
\drawpath{50.0}{8.0}{36.0}{8.0}
\drawpath{36.0}{8.0}{34.0}{12.0}
\drawpath{18.0}{26.0}{58.0}{26.0}
\drawpath{68.0}{26.0}{70.0}{26.0}
\drawpath{70.0}{26.0}{52.0}{8.0}
\drawpath{52.0}{8.0}{48.0}{4.0}
\drawpath{48.0}{4.0}{38.0}{4.0}
\drawpath{38.0}{4.0}{36.0}{4.0}
\drawpath{36.0}{4.0}{32.0}{14.0}
\drawpath{32.0}{14.0}{38.0}{18.0}
\drawpath{32.0}{16.0}{30.0}{18.0}
\drawpath{38.0}{20.0}{54.0}{20.0}
\drawpath{66.0}{20.0}{82.0}{20.0}
\drawpath{38.0}{18.0}{52.0}{18.0}
\drawpath{30.0}{18.0}{34.0}{20.0}
\drawpath{34.0}{22.0}{30.0}{22.0}
\drawpath{34.0}{22.0}{58.0}{22.0}
\drawpath{68.0}{22.0}{82.0}{22.0}
\drawpath{64.0}{18.0}{82.0}{18.0}
\drawpath{34.0}{20.0}{38.0}{20.0}
\drawpath{30.0}{22.0}{32.0}{20.0}
\drawframebox{87.0}{21.0}{10.0}{10.0}{$S$}
\drawcenteredtext{10.0}{26.0}{$R$}
\end{picture}
}
\end{center}

Naturality gives equality to

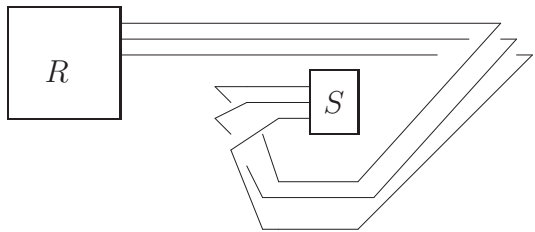
\begin{center}
{\tt\setlength{\unitlength}{3.0pt}
\begin{picture}(74,36)
\thinlines
\drawframebox{11.0}{25.0}{14.0}{14.0}{}
\drawpath{18.0}{30.0}{66.0}{30.0}
\drawpath{66.0}{30.0}{48.0}{10.0}
\drawpath{48.0}{10.0}{38.0}{10.0}
\drawpath{38.0}{10.0}{36.0}{16.0}
\drawpath{18.0}{28.0}{62.0}{28.0}
\drawpath{66.0}{28.0}{68.0}{28.0}
\drawpath{68.0}{28.0}{50.0}{8.0}
\drawpath{50.0}{8.0}{36.0}{8.0}
\drawpath{36.0}{8.0}{34.0}{12.0}
\drawpath{18.0}{26.0}{58.0}{26.0}
\drawpath{68.0}{26.0}{70.0}{26.0}
\drawpath{70.0}{26.0}{52.0}{8.0}
\drawpath{52.0}{8.0}{48.0}{4.0}
\drawpath{48.0}{4.0}{38.0}{4.0}
\drawpath{38.0}{4.0}{36.0}{4.0}
\drawpath{36.0}{4.0}{32.0}{14.0}
\drawpath{32.0}{14.0}{38.0}{18.0}
\drawpath{32.0}{16.0}{30.0}{18.0}
\drawpath{34.0}{22.0}{42.0}{22.0}
\drawpath{38.0}{20.0}{42.0}{20.0}
\drawpath{38.0}{18.0}{42.0}{18.0}
\drawpath{30.0}{18.0}{34.0}{20.0}
\drawpath{34.0}{22.0}{30.0}{22.0}
\drawframebox{45.0}{20.0}{6.0}{8.0}{$S$}
\drawpath{34.0}{20.0}{38.0}{20.0}
\drawpath{30.0}{22.0}{32.0}{20.0}
\drawcenteredtext{10.0}{24.0}{$R$}
\end{picture}
}
\end{center}
and then to 

\begin{center}
{\tt\setlength{\unitlength}{3.0pt}
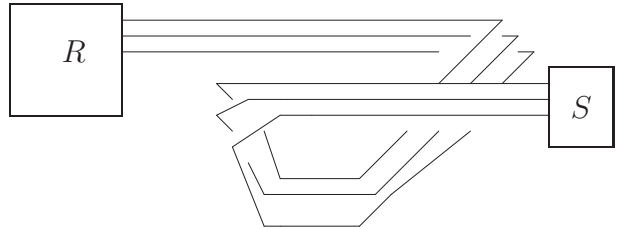
\begin{picture}(84,36)
\thinlines
\drawframebox{11.0}{25.0}{14.0}{14.0}{}
\drawpath{18.0}{30.0}{66.0}{30.0}
\drawpath{66.0}{30.0}{58.0}{22.0}
\drawpath{48.0}{10.0}{38.0}{10.0}
\drawpath{38.0}{10.0}{36.0}{16.0}
\drawpath{18.0}{28.0}{62.0}{28.0}
\drawpath{66.0}{28.0}{68.0}{28.0}
\drawpath{68.0}{28.0}{62.0}{22.0}
\drawpath{50.0}{8.0}{36.0}{8.0}
\drawpath{36.0}{8.0}{34.0}{12.0}
\drawpath{18.0}{26.0}{58.0}{26.0}
\drawpath{68.0}{26.0}{70.0}{26.0}
\drawpath{70.0}{26.0}{66.0}{22.0}
\drawpath{52.0}{8.0}{48.0}{4.0}
\drawpath{48.0}{4.0}{38.0}{4.0}
\drawpath{38.0}{4.0}{36.0}{4.0}
\drawpath{36.0}{4.0}{32.0}{14.0}
\drawpath{32.0}{14.0}{38.0}{18.0}
\drawpath{32.0}{16.0}{30.0}{18.0}
\drawpath{34.0}{22.0}{42.0}{22.0}
\drawpath{38.0}{20.0}{42.0}{20.0}
\drawpath{38.0}{18.0}{42.0}{18.0}
\drawpath{30.0}{18.0}{34.0}{20.0}
\drawpath{34.0}{22.0}{30.0}{22.0}
\drawpath{42.0}{22.0}{72.0}{22.0}
\drawpath{34.0}{20.0}{38.0}{20.0}
\drawpath{30.0}{22.0}{32.0}{20.0}
\drawcenteredtext{12.0}{26.0}{$R$}
\drawpath{42.0}{20.0}{72.0}{20.0}
\drawpath{42.0}{18.0}{72.0}{18.0}
\drawpath{62.0}{16.0}{52.0}{8.0}
\drawpath{58.0}{16.0}{50.0}{8.0}
\drawpath{54.0}{16.0}{48.0}{10.0}
\drawframebox{76.0}{19.0}{8.0}{10.0}{$S$}
\end{picture}
}
\end{center}

which is equal to the untwisted circuit.


\end{document}